\documentclass[11pt,twoside]{article}

\usepackage[a4paper,margin=3.5cm]{geometry}
\usepackage[T1]{fontenc}
\usepackage[utf8]{inputenc}
\usepackage[english]{babel}
\usepackage{stmaryrd} %for \llbracket and \rrbracket
\usepackage{csquotes}
\usepackage{booktabs}
\usepackage{amssymb}
\usepackage{amsmath}
\usepackage{amsthm}
\usepackage{dsfont}
\usepackage{enumitem}
\usepackage{epigraph}
\usepackage{fancyhdr}
\usepackage{mathrsfs}
\usepackage{mathtools}
\usepackage{graphicx}
\usepackage{array}
\usepackage{pgf,tikz}
\usepackage{tcolorbox}
\usepackage{tikz-cd}
\usepackage{titling}
\usepackage{url}
\usepackage{verbatim}
\usepackage[style=alphabetic, backend=biber]{biblatex}
\usepackage[colorlinks=true, allcolors=blue]{hyperref}
\usepackage{cleveref}

\setlist{noitemsep}

\usetikzlibrary{arrows}
\usetikzlibrary{decorations.pathmorphing}
\usetikzlibrary{patterns}

\newcommand{\cat}[1]{\mathcal #1}

\newcommand{\A}{\mathbb{A}}

\newcommand{\et}{\mathrm{\acute{e}t}}

\newcommand{\fppf}{\mathrm{fppf}}

\newcommand{\GG}{\mathbb{G}}

\newcommand{\hyptau}{\tn{hyp-}\tau}

\newcommand{\ideal}[1]{\mathfrak{#1}}

\newcommand{\lra}{\longrightarrow}

\newcommand{\mc}{\mathcal}

\newcommand{\mot}{\tn{mot}}

\newcommand{\N}{\mathbb{N}}
\newcommand{\Nis}{\mathrm{Nis}}

\renewcommand{\P}{\mathbb{P}}

\newcommand{\Q}{\mathbb{Q}}

\newcommand{\R}{\mathbb{R}}
\newcommand{\ra}{\rightarrow}

\let\parasign\S
\renewcommand{\S}{\mathbb{S}}

\newcommand{\stackstag}[1]{\cite[\href{https://stacks.math.columbia.edu/tag/#1}{Tag #1}]{stacks-project}}

\newcommand{\U}{\mathcal{U}}
\newcommand{\tn}{\text}
\newcommand{\transpose}[1]{{#1}^T}

\newcommand{\Z}{\mathbb{Z}}
\newcommand{\Zar}{\mathrm{Zar}}
\newcommand{\Gm}{\mathbb{G}_m}
\newcommand{\squotient}[2]{{#1/\!\!/#2}} %stacky quotient/homotopy quotient

\newcommand{\B}[1]{\mathrm{B}#1}

\newcommand{\characterlattice}{X^*}
\newcommand{\dominantcone}{X^+}
\newcommand{\cocharacterlattice}{X_*}
\newcommand{\coroot}[1]{{#1}^{\mathrm{co}}}
\DeclareMathOperator{\Aut}{Aut}
\DeclareMathOperator{\BG}{B{G}}

\DeclareMathOperator{\BO}{BO}
\DeclareMathOperator{\BSp}{BSp}

\DeclareMathOperator{\ch}{ch}
\DeclareMathOperator{\cha}{char}
\DeclareMathOperator*{\colim}{colim}

\DeclareMathOperator{\ESp}{ESp}

\DeclareMathOperator{\GL}{GL}

\DeclareMathOperator{\GW}{GW}

\DeclareMathOperator{\HGr}{HGr}

\DeclareMathOperator{\Hom}{Hom}
\DeclareMathOperator{\HP}{HP} % Quaternionic projective spaces
\DeclareMathOperator{\HU}{HU} % Universal torsors over quaternionic Grassmannians
\DeclareMathOperator{\id}{id}

\DeclareMathOperator{\K}{K}

\DeclareMathOperator{\res}{res}
\DeclareMathOperator{\rank}{rank}
\DeclareMathOperator{\sheafhom}{\mathscr{H}\textit{\kern -4pt om}\,}
\DeclareMathOperator{\Sing}{Sing}
\DeclareMathOperator{\SL}{SL}

\DeclareMathOperator{\Sp}{Sp}
\DeclareMathOperator{\Spec}{Spec}

\DeclareMathOperator{\Symp}{Symp}

\DeclareMathOperator{\Vect}{Vect}

\DeclareMathOperator{\Rep}{\mathbf{Rep}} % representation category
\DeclareMathOperator{\Sm}{\mathbf{Sm}}

\DeclareMathOperator{\sPre}{\mathbf{sPre}} % category of simplicial presheaves

\numberwithin{equation}{subsection}

\theoremstyle{definition}
\newtheorem{definition}{Definition}[subsection]

\theoremstyle{plain}
\newtheorem{corollary}[definition]{Corollary}
\newtheorem{lemma}[definition]{Lemma}
\newtheorem{proposition}[definition]{Proposition}
\newtheorem{theorem}[definition]{Theorem}

\theoremstyle{remark}
\newtheorem{example}[definition]{Example}

\newtheorem{remark}[definition]{Remark}

\title{\textsc{atiyah-segal completion for the hermitian k-theory of symplectic groups}}
\author{
    \textsc{jens hornbostel},
    \textsc{herman rohrbach},
    \textsc{marcus zibrowius}
}
\predate{}
\postdate{}
\date{}

\begin{document}

\maketitle

\begin{abstract}
    We study equivariant Hermitian K-theory for representations of symplectic groups, especially $\SL_2$. The results are used to establish an Atiyah-Segal completion theorem for hermitian $K$-theory and symplectic groups.
\end{abstract}
\tableofcontents

\section{Introduction}
    \label{section:introduction}

In \cite{rohrbach22completion}, a completion theorem for Hermitian K-theory of schemes with trivial torus action is established.
Let $X$ be a regular Noetherian separated scheme over $\Spec(\Z[\tfrac{1}{2}])$ with a trivial action of a split torus $T \cong \Gm^t$.
We let $IO = \ker \left(\GW^{T,[0]}_0(X) \ra \GW^{[0]}_0(X)\right)$ be the Hermitian version of the augmentation ideal.
Define
\begin{equation*}
    \GW^{T,[n]}(X)^{\wedge}_{IO}
\end{equation*}
as the derived completion of the $\GW^{T,[0]}(X)$-module $\GW^{T,[n]}(X)$ with respect to $IO$, see \cite{rohrbach22} for details.
Also define
\begin{equation*}
    \GW^{[n]}(\B{T}_X) = \lim_n \GW^{[n]}((\P^n_X)^{\times t}).
\end{equation*}
In \cite[Proposition 8.2.2]{rohrbach21}, it is shown that this definition is in fact an equivalence of $\GW$-spectra of motivic spaces.
In \cite[Theorem 3.2.3]{rohrbach22completion}, the following completion theorem for Hermitian K-theory is established.

\begin{theorem} \label{theorem:asc-for-split-tori}
    For all $i,n \in \Z$, the natural map
    \begin{equation*}
        \pi_i \left(\GW^{T,[n]}(X)^{\wedge}_{IO}\right) \stackrel{\cong}{\lra} \GW^{[n]}_i(\B{T}_X),
    \end{equation*}
    is an isomorphism.    
\end{theorem}

In the case of algebraic $\K$-theory, the corresponding result for split tori $T$, due to Totaro \cite{totaro99} for $\K_0$ and Knizel-Neshitov \cite{knizel14} for higher $\K$-groups, was the first important step for establishing more general versions of an Atiyah-Segal style completion theorem for linear algebraic groups, due to Krishna \cite{krishna18}, Tabuada-van den Bergh \cite{tabuada21} and Carlsson-Joshua \cite{carlsson23}, and their proofs all rely on a reduction to the case of a split torus.
The proof for algebraic K-theory and tori $T$ is easier than Rohrbach's theorem above on Hermitian $K$-theory because the classifying space $\B{T}$ is a product of copies of $\P^{\infty}$, and algebraic $K$-theory is orientable in the sense of Levine-Morel and Panin.
As of yet, there is no completion theorem for Hermitian $\K$-theory in the case of general linear groups $\GL_n$, although \cite{rohrbach23} contains a partial result in this direction by computing the Hermitian $\K$-theory of even-dimensional Grassmannians.

We restrict our attention to the Hermitian $\K$-theory of symmetric forms and symplectic forms in degree zero, which is the $\Z/2\Z$-graded Grothendieck-Witt group $\GW^{\pm}$.
To extend our results to higher Hermitian $\K$-theory over general base schemes, the second author intends to use the methods introduced in \cite{rohrbach22completion} in forthcoming work, using the machinery of derived completion, in line with the recent work of Tabuada-van den Bergh \cite{tabuada21} and Carlsson-Joshua \cite{carlsson23}.

It is likely difficult to extend the completion theorem for Hermitian $\K$-theory to schemes with non-trivial actions, as all the known proofs in algebraic $\K$-theory rely on the equivariant localization theorem, whose analogue in Hermitian $\K$-theory has not yet been established, and on a suitable geometric equivariant decomposition theorem. We refer to the last subsection for further details.

\medskip

In a series of fundamental articles, Panin and Walter establish the theory of symplectic oriented cohomology theories on smooth algebraic varieties. In this setting, the one-dimensional torus $\GL_1$ is replaced by $\Sp_2=\SL_2$. They also show that Hermitian K-theory is a symplectic orientable theory, which in particular implies the following computation of $\BSp_{2n}$ over a base field $k$ (see \cite{panin18} and \cite[Theorem 9.1]{panin22}):
\begin{equation*}
    \BO^{*,*}(\BSp_{2n}) \cong \BO^{*,*}(\Spec k)[[b_1, \dots, b_n]],
\end{equation*}
where the $b_i$ are the \emph{Borel classes} of \cite[Section 8]{panin22} and the right hand side is the ring of graded power series over $\BO^{*,*}(\Spec k)$, compare e.g.
\cite[section 6.3]{krishna12}.
From this point of view, one might argue that the conjectural Atiyah-Segal completion result for $\Sp_{2n}$ and Hermitian K-theory is the correct analogue of the completion theorem for $\GL_n$ and algebraic K-theory, and should provide computations involving free polynomial rings that are easier than Rohrbach's theorem above. 
We will show that this is indeed the case and prove the following theorem as Corollary \ref{corollary:asc-for-gw-of-sp} below.

\begin{theorem} \label{theorem:asc-for-gw-of-sp}
    Let $k$ be a field of characteristic not two.
    There is a canonical map of $\GW^{\pm}(k)$-algebras
    \begin{equation*}
        \GW^{\pm}(\Rep(\Sp_{2r})) \lra \GW^{\pm}(\BSp_{2r}),
    \end{equation*}
    which exhibits $\GW^{\pm}(\BSp_{2r})$ as the completion of $\GW^{\pm}(\Rep(\Sp_{2r}))$ with respect to $IO_{\Sp_{2r}}$. 
\end{theorem}

The strategy for $\GW$-theory consists in some sense in systematically replacing $\GL_1$ by $\Sp_2$, and more generally $\GL_n$ by $\Sp_{2n}$ in all steps of the proof for complex or algebraic $K$-theory. 

This approach has already been used in computations for Chow-Witt groups, see e.g. \cite{hornbostelwendt19}, and is motivated in part by the fact that real realization of $\Sp_{2n}$ has $U(n)$ as maximal compact subgroup. 

The techniques used in this proof, which build on and slightly generalize those of Morel-Voevodsky and Panin-Walter, can be used to construct geometric classifying spaces for other groups $G$ and other Grothendieck topologies $\tau$. 
A general result about geometric descriptions of classifying spaces $\B_{\tau}{G}$ in motivic homotopy theory is given in Theorem \ref{generalBGtheorem}.
 
\medskip
Readers familiar with the proof of \cite{atiyah69} might consider still another strategy for proving general Atiyah-Segal completion theorems for real and Hermitian $K$-theory, respectively: rather than $\Sp_{2n}$ and its subgroup $\Sp_2^{\times n}$, one should work with $\GL_n$ and its maximal split torus, but both equipped with a suitable involution that up to homotopy becomes complex conjugation under complex realization. 
There are indeed nice models for these in the algebro-geometric setting, which are explained in the introduction of Section \ref{section:gw-of-real-split-reductive-groups}.
In particular, the results in Section \ref{subsection:symmetric-representation-ring} apply to these algebraic groups with involution. 
However, we are not yet able to perform all necessary equivariant $\GW$-computations for these algebraic models of $\GL_n$ and split tori with involution, and the authors hope to return to this topic in the future.  

\medskip

This article, particularly \Cref{section:gw-of-real-split-reductive-groups}, contains results about the \emph{real} (sometimes \emph{quadratic} or \emph{Hermitian}) representation ring of a large class of reductive groups with involution, such as $\Sp_{2n}$ and $\GL_n$ with involution. This builds on previous work of Calm\`es-Hornbostel \cite{calmes05} and Zibrowius \cite{zibrowius15} and is of independent interest.  As an example of how some of the techniques of this paper extend to such real groups, we have included in \Cref{subsection:classifying-space-group-with-involution} a short construction of a classifying space for the multiplicative group \(\Gm\) with a non-trivial involution. 

\medskip

\emph{Acknowledgement}: This research was conducted in the framework of the DFG-
funded research training group GRK 2240: \emph{Algebro-Geometric Methods in Algebra,
Arithmetic and Topology}.
The second author was supported by the ERC through the project QUADAG. 
This paper is part of a project that has received funding from the European Research Council (ERC) under the European Union's Horizon 2020 research and innovation programme (grant agreement No. 832833).
The first author thanks Marc Hoyois for explaining several details in \cite{hoyois20cdh}.
The second author would like to thank Chris Schommer-Pries and Denis-Charles Cisinski for helpful comments on model structures on simplicial presheaves. 
We also kindly thank the anonymous referee for their detailed comments and suggestions.

\section{The Grothendieck-Witt rings of some Real split reductive groups}
\label{section:gw-of-real-split-reductive-groups}

Recall that the proof of Atiyah-Segal completion for topological KO-theory in \cite{atiyah69} relies on the splitting of \cite[Proposition~(5.2)]{atiyah68} for topological KR-theory for Real groups, i.e.\ for groups with an involution. We now discuss the algebraic analogue of such groups. 

A \emph{Real} group is a an algebraic group \(G\) with an involution \(\iota\colon G\to G\).  This involution is assumed to be a group homomorphism.  We denote by \(\Rep(G)\) the abelian category of finite-dimensional representations of \(G\), and equip it with a duality \(\vee_\iota\) as follows.
Recall first that, irrespective of any involution on \(G\), associating with a representation \(E\) of \(G\) the dual representation \( E^\vee := \Hom_G(E,k)\) defines a duality  \(\vee\) on \(\Rep(G)\).  More precisely, we obtain a category with duality \((\Rep(G), \vee, \eta)\), where \(\eta\) denotes the canonical double-dual identification.
Given the involution \(\iota\) on \(G\), we define the associated duality \(\Rep(G)\) as the composition \(\vee_\iota := \vee \circ \iota^*\), i.e.\ \( E^{\vee_{\iota}} := (\iota^* E)^\vee\).  We thus obtain a category with duality \(\Rep(G,\iota) := (\Rep(G), \vee_\iota,\eta_\iota)\).

In this section, we compute the Grothendieck-Witt ring
  \[
    \GW^\pm(\Rep(G,\iota)) := \GW^+(\Rep(G,\iota)) \oplus \GW^-(\Rep(G,\iota))
  \]
in cases when \(G\) is split reductive, the involution \(\iota\) restricts to a maximal torus \(T\), and all irreducible representations are self-dual with respect to \(\vee_\iota\).  

\newcommand{\inv}{\mathrm{inv}}
\begin{example}\label{eg:torusinv}
  Consider \((T,\inv)\), where \(T\) is a split torus of rank~\(r\) and \(\inv\) is the involution given by \(z\mapsto z^{-1}\). All representations of \(T\) symmetric with respect to~\(\vee_\inv\).
\end{example}

\begin{example}\label{eg:glninv}
  Consider \((\GL_n,\iota)\) with \(\iota(A) := \transpose{(A^{-1})}\), 
  an involution that restricts to the involution \(\inv\) on the standard maximal torus of \(\GL_n\).
  All representations of \(\GL_n\) are symmetric with respect to \(\vee_\iota\).
\end{example}

\begin{example}\label{eg:sp2n}
  Consider \((\Sp_{2n},\id)\).  All representations of \(\Sp_{2n}\) are self-dual with respect to \(\vee_{\id}\).  Some are symmetric, some are anti-symmetric (see \Cref{expl:sp} for more details).
\end{example}

\subsection{Reductive groups -- the setup}
\label{subsection:reductivegroups}
Let $k$ be a field. 
Let $(G, B, T)$ be a triple consisting of a connected split reductive group $G$ over $k$ with $T \subset B \subset G$, where $B$ is a Borel subgroup and $T$ is a maximal torus of rank $t$.
We write \(\characterlattice := \Hom(T, \GG_m)\) and \(\cocharacterlattice := \Hom(\GG_m, T)\) for the character lattice and the cocharacter lattice, respectively, both isomorphic to \(\Z^t\), and \(\langle -, - \rangle\colon \characterlattice \times \cocharacterlattice \to \Z\) for the canonical pairing between them.
Let $\Phi = \Phi(G,T)$ denote the associated set of roots, $\Phi^+ = \Phi(B,T)$ the set of positive roots associated with the Borel subgroup $B$ as in \cite[Proposition 1.4.4]{conrad14}, and $\Delta \subset \Phi^+$ the set of simple positive roots.  We write \(\coroot{\alpha}\) for the coroot associated with a root \(\alpha\).\footnote{
  We deviate from the universally agreed notation \(\alpha^\vee\) for \(\coroot{\alpha}\) to avoid any confusion with the duality \(\vee\) on \(\Rep(G)\).
}
Let $\dominantcone\subset \characterlattice$ denote the cone of dominant characters determined by \(\Phi^+\), also known as the (closed, integral) fundamental Weyl chamber \cite[Equation~(1.5.3)]{conrad14}.  Explicitly, the relation between \(\Phi^+\) and \(\dominantcone\) can be described as follows:
\begin{align}
  \dominantcone &= \{ x \in \characterlattice \mid \langle x, \coroot{\alpha} \rangle \geq 0 \text{ for all } \alpha \in \Delta \}\\
  \Phi^+        &= \{ \alpha \in \Phi \mid \langle x, \coroot{\alpha} \rangle \geq 0 \text{ for all } x \in \dominantcone \} \label{eq:positive-roots-in-terms-of-dominant-weights}
\end{align}

\begin{table}
  \begin{center}
    \begin{tabular}{cccl}
      \toprule
      \textit{here}
      & \cite{conrad14}
      & \cite{serre68}
      &
      \\
      \midrule
      $\characterlattice$
      & $X(T)$
      & $M$
      & character lattice
      \\
      $\dominantcone$
      & $C$
      & $P$
      & \(\substack{\text{dominant cone}\slash \\\text{[fundamental] Weyl chamber}}\)   \\
      $\Phi$
      & $\Phi$
      & $R$
      & set of roots
      \\
      $\Phi^+$
      & $\Phi^+$
      & $R^+$
      & set of positive roots
      \\
      $\Delta$
      & $\Delta$
      & ---
      & set of simple roots
      \\
      $\cocharacterlattice$
      & $X^\vee(T)$
      & ---
      & cocharacter lattice
      \\
      \bottomrule
    \end{tabular}
  \end{center}
  \caption{Translation between notation used here and notation in some of our references.}
  \label{translation}
\end{table}
The dominant characters parametrize the irreducible representations of $G$ \cite[lemme 5]{serre68}:  for each $x \in \dominantcone$, there is an irreducible \(G\)-representation $E_{x}$ of highest weight $x$, unique up to isomorphism.  We fix one such representation for each \(x\in \dominantcone\).
When \(G\) is simply connected, $\dominantcone \cong \N_0^n$, with basis given by fundamental weights $\omega_1, \dots, \omega_n$.  Following \cite[§\,3.6]{serre68}, we define a partial order on \(\characterlattice\) as follows:

\begin{definition}\label{def:serres-partial-order}
  For \(x,y\in \characterlattice\), we write $x \leq y$ if and only if $y - x$ can be written as a $\Z$-linear combination of elements of $\Phi^+$ with non-negative coefficients.
\end{definition}
(Note that other definitions abound in the literature.  The ordering defined here is slightly different from both orderings defined in \cite[Chapter~VI, Definition~2.2]{broeckertomdieck}.)

Finally, let $W = W(G,T) = W(\Phi)$ be the associated Weyl group.  We refer to any translate \(w\dominantcone \subset \characterlattice\) as a Weyl chamber. The Weyl group acts simply transitively on the set of Weyl chambers.  For the convenience of the reader, a partial translation of the notation used here and the notation used in some of the references is provided by \Cref{translation}.

\subsection{The representation ring}
\label{subsection{therepresentationring}}
\begin{lemma}\label{lemma:highestweightofproductissumofweights}
For dominant weights $x,y \in \dominantcone$,
\begin{equation*}
  E_x\cdot E_y = E_{x+y} + \textstyle\sum_z E_z 
\end{equation*}
in \(\K_0(\Rep(G))\), where the sum is over a finite number of \(z\in \dominantcone\) such that $z < x+y$.
\end{lemma}
\begin{proof}
  By \cite[th{\'e}or{\`e}me 4]{serre68}, there is an injective ring morphism 
  \begin{equation*}
      \ch_G: \K_0(\Rep(G)) \ra \Z[\characterlattice]
  \end{equation*}
  with image $\Z[\characterlattice]^W$. 
  By \cite[lemme 5]{serre68},
  \(
  \ch_G(E_x) = e^x + \sum_i e^{x_i}
  \)
  for a finite number of \(x_i \in \characterlattice\) such that $x_i < x$, and similarly for \(\ch_G(E_y)\) and \(\ch_G(E_{x+y})\).  It follows that 
\begin{equation*}
  \ch_G(E_x) \ch_G(E_y) = e^{x+y} + (\text{terms smaller than \(x+y\)}),
\end{equation*}
where by “terms smaller than \(x+y\)” we mean a \(\Z\)-linear combination of terms \(e^z\) with \(z\in \characterlattice\) such that \(z<x+y\).
On the other hand, \(\ch_G(E_x) \ch_G(E_y) = \ch_G(E_x \otimes E_y)\), and \(E_x \otimes E_y\) can be decomposed into a sum of irreducible representations \(E_1,\dots, E_N\).
Writing \(\ch_G(E_k) = e^{w_k} + (\text{terms smaller than} w_k)\), we find that  
\begin{equation*}
    \ch_G(E_x) \ch_G(E_y) = \textstyle\sum_{k = 1}^{N} \left( e^{w_k} + \text{terms smaller than \(w_k\)}\right)
\end{equation*}
It follows by comparison that \(w_k\leq x+y\) for each \(k\), and that \(w_k = x+y\) for exactly one \(k\). This proves the claim.
\end{proof}

\begin{lemma}\label{lemma:no-dominant-smaller-zero}
  There is no \(z\in \dominantcone\) such that \(z < 0\).
\end{lemma}
\begin{proof}
  Let \(\Delta\subset \Phi^+\) denote the set of simple roots.  
  If \(z<0\) then \(z = \sum_{\alpha\in \Delta}a_\alpha\alpha\) with \(a_\alpha \leq 0\) for all \(\alpha\in\Delta\) and \(a_\alpha < 0 \) for at least one \(\alpha\).
  In particular, \(z\) lies in the span of \(\Phi\), so we can pass from the root system \(\Phi\subset X^*\otimes_\Z \R\) to the \emph{reduced} root system \(\Phi \subset \Phi\otimes_\Z \R\), to which the results of \cite[chap.~VI, §\,1, n$^{\circ}$10]{bourbakiLie456} apply. 
  Let \(\sigma\) be the half-sum of all coroots \(\coroot{\alpha}\) corresponding to roots \(\alpha \in \Phi^+\).
  Then by  \cite[chap.~VI, §\,1, n$^{\circ}$10, corollaire]{bourbakiLie456}, \(\langle \alpha, \sigma \rangle = \sum_{\alpha\in\Delta}a_\alpha\).  This number is smaller than zero by assumption, so \(z\not\in \dominantcone\).
\end{proof}

\begin{proposition}\label{prop:representationringispolynomialring}
  Suppose \(\dominantcone\) can be split into a direct sum \(\dominantcone \cong \N_0^n \oplus \Z^m\).  Pick elements \(\omega_1,\dots,\omega_n\in \dominantcone\) and \(\zeta_1,\dots,\zeta_m\in \dominantcone\) corresponding to an \(\N_0\)-basis of \(\N_0^n\) and a \(\Z\)-basis of \(\Z^m\), respectively.  There is a unique ring homomorphism
  \[
    \Z[w_1,\dots,w_n,z_1^{\pm 1},\dots,z_m^{\pm 1}] \to \K_0(\Rep(G))
  \]
  taking each \(w_i\) to \(E_{\omega_i}\) and each \(z_j\) to \(E_{z_j}\), and this ring homomorphism is an isomorphism.
\end{proposition}
\begin{example}
  For simply connected \(G\), \(\dominantcone\cong \N_0^n\) with basis \(\omega_1,\dots,\omega_n\) the fundamental weights. The proposition shows that \(\K_0(\Rep(G))\) is a polynomial ring over \(\Z\) on generators \(E_{\omega_1},\dots,E_{\omega_n}\).
\end{example}
\begin{example}
  For \(G=T\) a split torus of rank \(m\), \(\dominantcone = \characterlattice \cong \Z^m\). Pick a \(\Z\)-basis \(\zeta_1,\dots,\zeta_m\) of \(\characterlattice\). The proposition shows that \(\K_0(\Rep(T))\) is a ring of Laurent polynomials with the one-dimensional representations \(E_{\zeta_i}\) as generators.
\end{example}
\begin{proof}[Proof of \Cref{prop:representationringispolynomialring}]
  It suffices to show that there is a well-defined ring map
  \[
    f\colon \frac{\Z[w_1,\dots,w_n,z_1,\dots,z_m,z_1',\dots,z_m']}{(z_jz_j'-1 \mid j = 1,\dots,m)} \to \K_0(\Rep(G))
  \] 
  that sends \(w_i\) to \(E_{\omega_i}\), \(z_j\) to \(E_{\zeta_j}\) and \(z_j'\to E_{-\zeta_j}\), and that this map is an isomorphism.

  To see that \(f\) is well-defined, note that by \Cref{lemma:highestweightofproductissumofweights}
  \[
    f(z_j)f(z_j')
    = E_{\zeta_j}\cdot E_{-\zeta_j}
    = E_0 + \textstyle\sum_z E_z,
  \]
  where the sum is over certain \(z\in \dominantcone\) such that \(z < 0\).  By \Cref{lemma:no-dominant-smaller-zero}, this is the empty sum.  So \(f(z_j)f(z_j') = E_0 = 1\), as required.
  
  To see that \(f\) is an isomorphism, we associate the Laurent monomial \(M^x := w_1^{a_1} \cdot \dots \cdot w_n^{a_n}\cdot z_1^{b_1}\cdot \cdots \cdot z_m^{b_m}\) with the element \(x = \sum_i a_i \omega_i + \sum_j b_j \zeta_j \in \dominantcone\). This defines a bijection between a \(\Z\)-basis of the Laurent ring and \(\dominantcone\). Under \(f\), the basis element \(M^x\) maps to
  \(
  \prod_i(E_{\omega_i})^{a_i}\prod_j(E_{\zeta_j})^{b_j}
  \).
  By \Cref{lemma:highestweightofproductissumofweights}, we can rewrite this element as 
  \begin{equation}\label{eq:representationringispolynomialring:key-identity}
    f(M^x) = E_x + \textstyle\sum_z E_z,
  \end{equation}
  for certain \(z\in \dominantcone\) with \(z<x\).  
  Using \cite[chap.~VI, §\,3, n$^{\circ}$4, lemme~4]{bourbakiLie456} as in the proof of \cite[Lemme~6]{serre68},
  we deduce that the elements \(f(M^x)\) for \(x\in \dominantcone\) form a \(\Z\)-basis of \(\Z[\characterlattice]^W\).
  This concludes the proof.
\end{proof}

\begin{remark}[Alternative generators]\label{rem:representationringispolynomialring-alternative-generators}
  More generally, under the assumptions of \Cref{prop:representationringispolynomialring}, we can choose generators for \(\K_0(\Rep(G))\) as follows.
  Take \(\omega_i\) and \(\zeta_j \in \dominantcone\) as before.  For \(i\in \{1,\dots,n\}\), pick classes \(e_{\omega_i} \in \K_0(\Rep(G))\) such that 
  \[
    e_{\omega_i} = E_{\omega_i} + (\text{smaller terms}),
  \]
  in \(\K_0(\Rep(G))\), where (smaller terms) refers to a \(\Z\)-linear combination of irreducible representations \(E_x\) indexed by finitely many \(x \in \dominantcone\) with \(x < \omega_i\).
  Then again we have a ring isomorphism 
  \[
    \Z[w_1,\dots,w_n,z_1^{\pm 1},\dots,z_m^{\pm 1}] \xrightarrow{\cong} \K_0(\Rep(G))
  \]
  taking each \(w_i\) to \(e_{\omega_i}\) and each \(z_j\) to \(E_{\zeta_j}\).
  Indeed, this follows with the same proof, as the key identity \eqref{eq:representationringispolynomialring:key-identity} still holds for these more general generators.
\end{remark}

\subsection{The \texorpdfstring{$\pm$}{±}-symmetric representation ring}
    \label{subsection:symmetric-representation-ring}

From now on, we always assume $\cha(k) \neq 2$.

We now assume that \(G\) is equipped with an involution \(\iota\) (possibly trivial) \emph{that restricts to the chosen maximal torus \(T\)}, and study the associated duality \(\vee_\iota\) on \(\Rep(G)\) introduced at the beginning of the section.

The following lemmas describe the highest weight of the \(\vee_\iota\)-dual of an irreducible representation.  As we have assumed that \(\iota\) restricts to \(T\), we have induced involutions \(\iota^*\) and \(\iota_*\) on the character lattice \(\characterlattice\) and the cocharacter lattice \(\cocharacterlattice\), respectively, compatible with the canonical pairing \(\langle -, -\rangle\).  For example, the involutions \(\inv^*\) and \(\inv_*\) induced by the involution from \Cref{eg:torusinv} are given by \(-\id\).  Interpreting characters as one-dimensional representations, we see that the dualities \(\vee\) and \(\vee_\iota\) also define involutions on \(\characterlattice\).  Explicitly, \(\vee = -\id\) on \(\characterlattice\), hence \(\vee_\iota = -\iota^*\) on \(\characterlattice\).

\begin{lemma}\label{lemma:Weyl-element-enconding-iota-duality}
  There is a unique involution \(w_\iota\in W\) whose action on Weyl chambers agrees with this action of \(\vee_\iota = -\iota^*\) on the Weyl chambers.
  In particular, \(-\omega_\iota\iota^*\) defines an involution of the fundamental Weyl chamber \(\dominantcone\).
\end{lemma}
\begin{example}
  If \(\iota|_T = \id\), then \(\vee_\iota = -\id\) on \(\characterlattice\), and \(w_\iota = w_0\), the longest element of the Weyl group.
  If \(\iota|_T = \inv\) is the inversion from \Cref{eg:torusinv}, then \(\vee_\iota = \id\) on \(\characterlattice\), and \(w_\iota = 1\).
\end{example}
\begin{proof}
  The case when \(\iota = \id\) is standard.
  For the general case, note that both involutions \(\iota^*\) and \(\vee\) on \(\characterlattice\) send roots to roots.
  (For \(\iota^*\), this follows from \(\iota\) is a group homomorphism compatible with the inclusion \(T\subset G\); c.f.\ \cite[above Theorem~1.3.15]{conrad14}.)
  The claim then follows from the fact that the Weyl group acts simply transivitely on the set of Weyl chambers.
\end{proof}

\begin{lemma}\label{lemma:weightofdualrepresentation}
  Given an irreducible $G$-representation \(E_x\) with highest weight \(x\), the dual representation \((\iota^*E_x)^\vee\) has highest weight \(-w_\iota \iota^*x\), where \(w_\iota\) is as in the previous lemma.
\end{lemma}
\begin{proof}
  Let $\Omega_E \subset \characterlattice$ be the set of weights of $E := E_x$.
  As \(\vee_\iota = -\iota^*\) on \(\characterlattice\), we have $\Omega_{(\iota^*E)^\vee} = -\iota^*\Omega_{E}$.
  Now consider an arbitrary weight \(z\in\Omega_{(\iota^*E)^\vee}\).
  As \(\Omega_{(\iota^*E)^\vee}\) is invariant under the action of the Weyl group, we also have \(\omega_\iota z \in \Omega_{(\iota^*E)^\vee}\).  
  So \(\omega_\iota z \in (-\iota^*\Omega_{E})\).  
  As $x$ is the highest weight of $E$, this means that 
  \(
    \omega_\iota z = -\iota^*(x - \sum_{b \in \Delta} m_b b)
  \)
  for certain \(m_b \in \Z_{\geq 0}\).  So
  \(
  z = -\omega_\iota\iota^*x + \sum_{b \in w_\iota\iota_*(\Delta)} m'_b b
  \)
  for certain \(m'_b\in\Z_{\geq 0}\).
  As noted in \Cref{lemma:Weyl-element-enconding-iota-duality}, \(\omega_\iota\iota^*(\dominantcone) = -\dominantcone\).
  In view of \Cref{eq:positive-roots-in-terms-of-dominant-weights}, this implies that \(w_\iota\iota_*(\Delta)\subseteq -\Phi^+\).  
  So
  \(
  z = -\omega_\iota\iota^*x - \sum_{b\in\Phi^+} m''_bb
  \),
  for certain \(m''_b \in\Z_{\geq 0}\).  It follows that $-w_\iota\iota^*x$ is the highest weight of \((\iota^*E)^\vee\), as claimed.
\end{proof}

\begin{definition} \label{definition:self-dualcharacter}
  A dominant character $x \in \dominantcone$ is (\(\vee_\iota\)-)\emph{self-dual} if $x = -w_\iota\iota^* x$, where $w_\iota$ is as in \Cref{lemma:Weyl-element-enconding-iota-duality}.
\end{definition}

\begin{definition} \label{definition:self-dualrepresentation}
We call a representation $E$ of $G$ \emph{self-dual} if there exists an isomorphism $\phi\colon E \ra E^{\vee_\iota}$. 
We call a self-dual representation $E$ \emph{symmetric} or \emph{anti-symmetric} if $\phi$ can be chosen to be symmetric or anti-symmetric, respectively.
\end{definition}
By \Cref{lemma:weightofdualrepresentation}, the irreducible representation with highest weight \(x\) is \(\vee_\iota\)-self-dual if and only if \(x\) is \(\vee_\iota\)-self-dual. 
In general, a representation may be both symmetric and anti-symmetric.
However, in our setting a self-dual \emph{irreducible} representation is always either symmetric or anti-symmetric, but never both:

\begin{lemma} \label{lemma:self-dualsimpleissymmetricoranti-symmetric}
Let $E$ be an irreducible representation of a split reductive group $G$. 
If there exists an isomorphism $\phi\colon E \ra (\iota^*E)^\vee$, then $\phi$ is either symmetric or anti-symmetric, and any isomorphism $E \ra (\iota^*E)^\vee$ is a multiple of $\phi$ by an (invertible) scalar.
\end{lemma}
\begin{proof}
  This is essentially \cite[corollary 2.3]{calmes05}, which relies on \cite[Lemma 1.21]{calmes05}.
  Note that \cite[§1]{calmes05} is phrased in the generality of abelian categories with duality and hence applies to the category of \(G\)-representations regardless of our choice of duality on this category.
\end{proof}

In light of this lemma, we make the following definition.

\begin{definition} \label{definittion:signofcharacter}
The \emph{sign} \(s(x)\) of a dominant character $x \in \dominantcone$ is an element of $\{-1,0,1\}$, which is $1$ if $E_x$ is symmetric, $-1$ if $E_x$ is anti-symmetric, and $0$ if $E_x$ is not \(\vee_\iota\)-self-dual.
\end{definition}

For an abelian category \(\cat A\) with duality, i.e.\ with a fixed involution \(\vee\) and double-dual identification \(\omega\), we write \(\GW^+(\cat A,\vee,\omega)\) and \(\GW^-(\cat A,\vee,\omega)\) for the Grothendieck-Witt groups of symmetric and anti-symmetric forms over \(\cat A\), respectively.  Taking \(\cat A\) to be the category of finite-dimensional \(k\)-vector spaces with its usual duality, this construction yields the usual groups \(\GW^+(k)\) and \(\GW^-(k)\) of symmetric and anti-symmetric bilinear forms over \(k\).  The tensor product furnishes us with a \(\Z/2\)-graded ring structure on the direct sum \(\GW^\pm(k) := \GW^+(k) \oplus \GW^-(k)\).  As it is convenient to think of the grading group \(\Z/2\) multiplicatively as \(\{\pm 1\}\), we will refer to this grading as a \(\pm\)-grading in all that follows.
Note that \(\GW^-(k)\) does not depend on the field \(k\). Explicitly, \(\GW^-(k) = \Z\cdot H^-\) as a group, and 
\begin{equation}\label{eq:graded-GW-of-field}
  \GW^\pm(k) = \frac{\GW^+(k)[H^-]}{((H^-)^2 - 2H^+)}
\end{equation}
as a ring, where \(H^\pm\in\GW^\pm(k)\) are the respective hyperbolic planes.

More generally, for the category with duality \((\Rep(G),\iota)\) introduced at the beginning of this section, we obtain a \(\pm\)-graded algebra \(\GW^\pm(\Rep G,\iota)\) over the \(\pm\)-graded ring \(\GW^\pm(k)\). 
The following theorem is an analogue of \cite[Theorem~2.10]{calmes05} for Grothendieck-Witt theory, under the assumption that all characters are self-dual so as to eliminate the existence of hyperbolic elements.

\begin{proposition}\label{prop:generatorsofgradedgwring}
Assume that all $x \in \dominantcone$ are \(\vee_\iota\)-self-dual.
Choose an isomorphism $\phi_x\colon E_x \ra (\iota^*E_x)^\vee$ for each $x \in \dominantcone$. 
Then the classes $[E_x, \phi_x]$ form a basis of the $\GW^{\pm}(k)$-module $\GW^{\pm}(\Rep(G,\iota))$. 
\end{proposition}
\begin{proof}
  The proof is the same as for \cite[Theorem~2.10]{calmes05}, noting that the obstruction from \cite[Remark~2.11]{calmes05} has been removed with the assumption that all dominant characters are self-dual.
  The main ingredients are \cite[Corollary~1.14]{calmes05}, which provides additive decompositions of both \(\GW^+(\Rep(G))\) and \(\GW^-(\Rep(G))\), and \cite[Corollary~1.38]{calmes05}, which identifies the summands of \(\GW^+(\Rep(G)\) corresponding to symmetric characters and the summands of \(\GW^-(\Rep(G))\) corresponding to anti-symmetric characters. For a full proof, \cite[Corollary~1.38]{calmes05} and the preceding \cite[Proposition~1.37]{calmes05} need to be mildly generalized to include all signs.  (For example, strictly speaking the identification of \(\GW(\mathcal A_i)\) with \(\GW_-(\mathcal A_{\mathbf{1}})\) as a \(\GW(\mathcal A_{\mathbf{1}})\)-module is missing from \cite[Corollary~1.38\,(i)]{calmes05}, as the tensor unit \(\mathbf{1}\) is assumed to be symmetric, not just \(\delta\)-symmetric, throughout \cite[\parasign\,1.4]{calmes05}.)  
\end{proof}

\begin{remark}
  Even in the presence of non-self-dual characters, it is easy to describe the \(\GW^\pm(k)\)-module structure of \(\GW^\pm(\Rep(G,\iota))\).
  Let \(H^+(E)\) and \(H^-(E)\) denote the symmetric and anti-symmetric hyperbolic space associated with a representation \(E\), respectively.
  The results quoted from \cite{calmes05} in the proof above show that every pair of non-self-dual dominant characters \((z,-w_\iota z)\) contributes a copy of \(\Z\) generated by \(H^+(E_z)\) to \(\GW^+(\Rep(G,\iota))\) and a copy of \(\Z\) generated by \(H^-(E_z)\) to \(\GW^-(\Rep(G,\iota))\).
  However, we want to concentrate on the case when \(\GW^\pm(\Rep(G,\iota))\) is a \emph{free} \(\GW^\pm(k)\)-module here.
\end{remark}

\begin{lemma}\label{lem:symmetric-rep-decomposition}
  Consider the duality \(\vee_\iota\) on \(\Rep{G}\).
  Suppose \((V,\phi)\) is an \(\epsilon\)-symmetric representation of \(G\) such that \(V = E_x + (\text{smaller terms})\) in \(\K_0(\Rep(G))\), for some \(x\in \dominantcone\) and \(\epsilon \in\{\pm\}\). 
  Then there exists an \(\epsilon\)-symmetric isomorphism \(\phi_x\) on \(E_x\) such that 
  \[
    (V,\phi) = (E_x,\phi_x) + (\text{smaller terms})
  \]
  in \(\GW^\epsilon(\Rep(G,\iota))\), where (smaller terms) refers to a \(\GW(k)\)-linear combination
  of terms \(T_z\) indexed by finitely many \(z \in \dominantcone\) with \(z < x\) such that, for each \(z\), \(T_z\) is either of the form \((E_z,\phi_z)\) (in case \(z\) is \(\epsilon\)-symmetric)
  or of the form \(H^{\epsilon}(E_z)\) (in case \(z\) is not).
\end{lemma}
\begin{proof}
  Using the decomposition of \(\GW^\epsilon(\Rep(G,\iota))\) of \cite[Corollary~1.14]{calmes05} and \cite[Corollary~1.38]{calmes05} as in the proof of \Cref{prop:generatorsofgradedgwring} above, 
  we can write \((V,\phi)\) as a \(\GW(k)\)-linear combination of elements \((E_z,\phi_z)\) with \(z\in \dominantcone\) such that \(\phi_z\) is \(\epsilon\)-symmetric, and of elements \(H^\epsilon(E_v)\) with \(v\in \dominantcone\) such that \(E_v\) is not \(\epsilon\)-self-dual:
  \[
    (V,\phi) = \sum_{\mathclap{\substack{z\colon s(z) = \epsilon}}} \alpha_z \cdot (E_z,\phi_z)
    + \sum_{\mathclap{\substack{v\colon s(v) = -\epsilon\\\text{ or } s(v) = 0}}} a_v\cdot H^\epsilon(E_v)
  \]
  The coefficients \(a_v\) may be taken in \(\Z\), as \(a \cdot H^{\epsilon}(b) =H^{\delta\epsilon}(F(a)b)\) for any \(a\in\GW^\delta(\Rep(G,\iota))\) and any \(b\in\K_0(\Rep(G))\).
  By \cite[Remark~1.15]{calmes05}, we even know that the coefficients \(\alpha_z\) of \((E_z,\phi_z)\) can be chosen to be actual symmetric forms over \(k\) of positive rank, and that the coeffients \(a_v\in\Z\) of \(H^\epsilon(E_v)\) are positive.
  Applying the forgetful map, we thus find that  
  \[
    V = \textstyle\sum_z (\rank{\alpha_z})E_z + \textstyle\sum_v 2a_v E_v
  \]
  in \(\K_0(\Rep(G))\).
  As the irreducible representations \(E_x\) with \(x\) ranging over \(\dominantcone\) form a \(\Z\)-basis of \(\K_0(\Rep(G))\), we deduce by comparing this decomposition with the given decomposition of \((V,\phi)\) that one of the \(z\)'s in the first sum must be equal to \(x\), with \(\alpha_x\) of rank one, and that for all other \(z\)'s and all \(v\)'s we have \(z < x\) and \(v < x\), respectively.
\end{proof}

The following is a careful restatement of \cite[Lemma~2.14]{calmes05}:
\begin{lemma}\label{lemma:calmes2.14}
  For any two self-dual \(x,y\in \dominantcone\), the sum \(x+y \in \dominantcone\) is also self-dual.
  Moreover, given \(\epsilon_x\)- and \(\epsilon_y\)-symmetric isomorphisms \(\phi_x\colon E_x\to (\iota^* E_x)^\vee\) and \(\phi_y\colon E_y \to (\iota^*E_y)^\vee\), respectively, there exists an \(\epsilon_x\epsilon_y\)-symmetric isomorphism \(\phi\colon E_{x+y}\to (\iota^* E_{x+y})^\vee\) such that
  \[
    (E_x,\phi_x)\cdot(E_y,\phi_y) = (E_{x+y},\phi) + (\text{smaller terms}),
  \]
  in \(\GW^{\epsilon_x\epsilon_y}(\Rep(G,\iota))\), where (smaller terms) is to be read as in \Cref{lem:symmetric-rep-decomposition}.
\end{lemma}
\begin{proof}
  The claim that \(x+y\) is self-dual follows immediately from the definitions.
  Also, \((E_x,\phi_x)\otimes(E_y,\phi_y)\) is an \(\epsilon_x\epsilon_y\)-symmetric representation and hence defines an element of \(\GW^{\epsilon_x\epsilon_y}(\Rep(G,\iota))\).
  From \Cref{lemma:highestweightofproductissumofweights}, we know that
  \(
  E_x \cdot E_y = E_{x+y} + (\text{smaller terms})
  \),
  so we can apply \Cref{lem:symmetric-rep-decomposition}.
\end{proof}

The following result should be compared to \cite[Theorem~2.16]{calmes05}.  It is an analogue, with more restrictive hypotheses, of the description of the usual representation ring provided in \Cref{prop:representationringispolynomialring}.

\begin{proposition}\label{prop:gwringispolynomialring}
  Suppose that all \(x\in \dominantcone\) are \(\vee_\iota\)-self-dual.  Assume in addition that \(\dominantcone\) can be split into a direct sum \(\dominantcone \cong \N_0^n \oplus \Z^m\).  Pick elements \(\omega_1,\dots,\omega_n\in \dominantcone\) and \(\zeta_1,\dots,\zeta_m\in \dominantcone\) corresponding to an \(\N_0\)-basis of \(\N_0^n\) and a \(\Z\)-basis of \(\Z^m\), respectively.
  Pick (anti-)symmetric isomorphisms \(\phi_{\omega_i}\) and \(\phi_{\zeta_j}\) for each of the representations \(E_{\omega_i}\) and \(E_{\zeta_j}\), respectively.
  Consider the ring
  \(
    \GW^{\pm}(k)[w_1,\dots,w_n,z_1^{\pm 1},\dots,z_m^{\pm 1}]
  \)
  as a graded \(\GW^\pm(k)\)-algebra with generators \(w_i\) and \(z_j\) in degrees \(s(w_i)\) and \(s(z_j)\), respectively.
  Then there is a unique graded \(\GW^{\pm}(k)\)-algebra homomorphism
  \[
    f: \GW^{\pm}(k)[w_1,\dots,w_n,z_1^{\pm 1},\dots,z_m^{\pm 1}] \to \GW^\pm(\Rep(G,\iota))
  \]
  taking each \(w_i\) to \((E_{\omega_i},\phi_{\omega_i})\) and each \(z_j\) to \((E_{\zeta_j},\phi_{\zeta_j})\), and this homomorphism is an isomorphism.
\end{proposition}
    Note that for the one-dimensional representations \(E_{\zeta_j}\) we have \(E_{\zeta_j}^\vee \cong E_{-\zeta_j}\) with respect to the usual duality \(\vee = \vee_{\id}\).  So, for \(m>0\), the assumption that each \(E_{\zeta}\) be self-dual with respect to \(\vee_\iota\) cannot hold with respect to the usual duality \(\vee\).

\begin{example}
  Take \((G,\iota)=(T,\inv)\), a split torus of rank \(m\) equipped with the inversion from \Cref{eg:torusinv}. In this case, \(\dominantcone = \characterlattice \cong \Z^m\), and all characters are self-dual symmetric. The proposition shows that \(\GW^\pm(\Rep(T),\inv)\) is a ring of Laurent polynomials over \(\GW^\pm(k)\) with one-dimensional symmetric representations as generators.
\end{example}

\begin{example}
  Take \((G,\iota)=(\GL_n,\iota)\), the general linear group equipped with the involution from \Cref{eg:glninv}.  The fundamental representations of \(\GL_n\) are given by the exterior powers of the standard \(n\)-dimensional representation \(V\) \cite[Part~II, \parasign\,2.15]{jantzen03representations}. More precisely, in the notation from above, \(E_{\omega_i} = \Lambda^i(V)\) for \(i = 1,\dots, n-1\) and \(E_{\zeta_1} = \Lambda^n(V)\). As noted in \Cref{eg:glninv}, we can choose a non-degenerate \(\vee_\iota\)-symmetric form \(\phi\) on \(V\), so that we obtain elements \(\lambda_i := (\Lambda^i(V),\Lambda^i(\phi)\) in \(\GW^+(\Rep(\GL_n,\iota))\).
  The above proposition then shows that
  \[
    \GW^\pm(\Rep(\GL_n,\iota)) \cong \GW^\pm(k)[\lambda_1,\dots,\lambda_{n-1},\lambda_n,\lambda_n^{-1}].
  \]
\end{example}

\begin{example}\label{expl:sp}
  Take \((G,\iota) = (\Sp_{2n},\id)\), as in \Cref{eg:sp2n}.
  The irreducible representation \(E_{\omega_i}\) is symmetric for even \(i\) and anti-symmetric for odd \(i\) \cite[vol.\,VIII,Table~1]{bourbakiLie789}, so we can choose a nondegenerate \((-1)^i\)-symmetric form \(\phi_i\) on each \(E_{\omega_i}\).
  The proposition then shows that \(\GW^\pm(\Rep(\Sp_{2n}),\id)\) is a polynomial algebra over \(\GW^\pm(k)\) on \(n\) generators \((E_{\omega_1},\phi_{\omega_1})\), \dots, \((E_{\omega_n},\phi_n)\), with \((E_{\omega_i},\phi_i)\) of degree \((-1)^i\). 
  (We will see in \Cref{lem:GW-Sp-alternative-generators}, using the generalization of \Cref{prop:gwringispolynomialring} discussed in \Cref{rem:gwringispolynomialring-alternative-generators}, that alternative polynomial generators are again given by the exterior powers of the standard representation.)
\end{example}

\begin{proof}[Proof of \Cref{prop:gwringispolynomialring}]
  We have seen in the proof of \Cref{prop:representationringispolynomialring} that \(E_{\zeta_j} \cdot E_{-\zeta_j}  = 1\) in \(\K_0(\Rep(G)\).
  In particular, each \(E_{\zeta_j}\) is one-dimensional.  We can therefore pick symmetric isomorphisms \(\phi_{-\zeta_j}\) on each \(E_{-\zeta_j}^\vee\) such that
  \[
    (E_{\zeta_j},\phi_{\zeta_j})\cdot(E_{-\zeta_j},\phi_{-\zeta_j}) = 1
  \]
  in \(\GW(\Rep(G,\iota))\).  It follows that we have a well-defined \(\GW^{\pm}(k)\)-algebra homomorphism
  \[
    f\colon \GW^{\pm}(k)[w_1,\dots,w_n,z_1^{\pm 1},\dots,z_m^{\pm 1}] \to \GW^\pm(\Rep(G,\iota))
  \]
  sending \(w_i\) to \((E_{\omega_i},\phi_{\omega_i})\), \(z_j\) to \((E_{z_j},\phi_{z_j})\), and \(z_j'\) to \((E_{-z_j},\phi_{-z_j})\).

  To see that this map is an isomorphism, we argue exactly as in the proof of \Cref{prop:representationringispolynomialring} and associate the Laurent monomial \(M^x := w_1^{a_1} \cdot \dots \cdot w_n^{a_n}\cdot z_1^{b_1}\cdot \cdots \cdot z_m^{b_m}\) with the element \(x = \sum_i a_i \omega_i + \sum_j b_j \zeta_j \in \dominantcone\).
  Under the map \(f\) above, the basis element \(M^x\) maps to
  \[
    \prod_{i = 1}^n (E_{\omega_i}, \phi_{\omega_i})^{a_i}\cdot \prod_{j=1}^m (E_{\zeta_j}, \phi_{\zeta_j})^{b_j},
  \]
  and it suffices to show that these classes form a basis of $\GW^{\pm}(\Rep(G))$ considered as $\GW^{\pm}(k)$-module.
  By \Cref{lemma:calmes2.14}, we can rewrite \(f(M^x)\) as
  \begin{equation}\label{eq:gwringispolynomialring:key-identity}
    f(M^x) = (E_x,\phi_x) + \sum_z a_z (E_z,\phi_z),
  \end{equation}
  where \(\phi_x\) is an appropriately chosen (anti-)symmetric isomorphism on \(E_x\), the sum is over certain \(z\in \dominantcone\) such that \(z < x\), and \(a_z \in\GW^\pm(k)\).  (Remember all \(z\in \dominantcone\) are self-dual by assumption.) \Cref{prop:generatorsofgradedgwring} tells us that, if we equip each irreducible representation \(E_x\) with the isomorphism \(\phi_x\) that appears here, the elements $(E_x, \phi_x)$ form a $\GW^{\pm}(k)$-module basis for $\GW^{\pm}(\Rep(G))$.
  The claim then again follows from \cite[chap.~VI, §\,3, n$^{\circ}$4, lemme~4]{bourbakiLie456}.
\end{proof}

Just as in the case of the usual representation ring (Remark~\ref{rem:representationringispolynomialring-alternative-generators}), we can and sometimes will pick more general generators than the ones specified in \Cref{prop:gwringispolynomialring}.

\begin{remark}[Alternative generators]\label{rem:gwringispolynomialring-alternative-generators}
  Under the assumptions of \Cref{prop:gwringispolynomialring}, we can more generally choose generators for \(\GW(\Rep(G,\iota))\) as follows.
  As before, pick (anti-)symmetric isomorphisms \(\phi_{\omega_i}\) and \(\phi_{\zeta_j}\) for each of the representations \(E_{\omega_i}\) and \(E_{\zeta_j}\), respectively  
  Pick homogeneous classes \(e_{\omega_i} \in \GW^\pm(\Rep(G,\iota))\) for \(i\in\{1,\dots,n\}\) such that
  \[
    e_{\omega_i} = (E_{\omega_i},\phi_{\omega_i}) + (\text{smaller terms}),
  \]
  where (smaller terms) is to be read as in \Cref{lem:symmetric-rep-decomposition}.
  Then again we have an isomorphism of graded \(\GW^{\pm}(k)\)-algebras
  \[
    f: \GW^{\pm}(k)[w_1,\dots,w_n,z_1^{\pm 1},\dots,z_m^{\pm 1}] \to \GW^\pm(\Rep(G,\iota))
  \]
  taking each \(w_i\) to \(e_{\omega_i}\) and each \(z_j\) to \((E_{\zeta_j},\phi_{\zeta_j})\).
  Indeed, this follows with the same proof, as the key identity \eqref{eq:gwringispolynomialring:key-identity} still holds for these more general generators.
\end{remark}

\subsection{Some details on symplectic groups}
\label{subsection:some-details-symplectic-groups}

Now let $G = \Sp_{2n}$ be the symplectic group, i.e.\ the group of automorphisms of the anti-symmetric form
\begin{equation}\label{eq:model-of-Sp}
  J := \begin{pmatrix}0 & \id_n \\ -\id_n & 0 \end{pmatrix},
\end{equation}
and \(T\) the maximal torus, of rank $n$, consisting of matrices
\begin{equation}\label{eq:T-of-Sp}
  \begin{pmatrix}
    D & 0 \\ 0 & D^{-1}
  \end{pmatrix}
\end{equation}
with $D \in \GL_n$ diagonal, see e.g.\ \cite[Exercise 2.4.6]{conrad14}.
The Weyl group of $G$ is the wreath product $S_2 \wr S_n$, that is, the semi-direct product $S_2^n \rtimes S_n$ with the permutation action of $S_n$ on $S_2^n$.
The character lattice $\characterlattice$ of $G$ is a free $\Z$-module of rank $n$, with dominant characters $\dominantcone \cong \N_0^n$. 

If $n = 1$, then $G = \SL_2=\Sp_2$ and $\dominantcone \cong \N_0$.
The standard 2-dimensional representation $V$ of $G$ preserves the standard nondegenerate anti-symmetric bilinear form $\phi$, yielding an element $[V, \phi] \in \GW^{-}(\Rep(G))$.
Note that $V$ is a simple representation of highest weight $1$, which is fundamental.
By \Cref{prop:gwringispolynomialring}, \(\GW^{\pm}(\Rep(\Sp_2))\) therefore is a polyonmial algebra over $\GW^{\pm}(k)$ on one generator \((V,\phi)\).
It will be convenient for us to use the rank zero element \(b := (V,\phi)-H^-\) as a generator instead, so that we have 
\begin{equation}\label{eq:GW-of-Sp2}
  \GW^{\pm}(\Rep(\Sp_2)) \cong \GW^{\pm}(k)[b]. 
\end{equation}
We think of \(b\) as a “first Borel class”. Note that \(b\) generates the graded augmentation ideal of \(\GW^{\pm}(\Rep(\Sp_2))\) (see \Cref{def:graded-augmentation-ideal}).

For \(n>1\), it will be convenient for us to replace the standard polynomial generators \(E_{\omega_i}\) of \(\K_0(\Rep(\Sp_{2n}))\) and the corresponding generators of \(\GW^\pm(\Rep(\Sp_{2n}))\) from \Cref{expl:sp} by alternative polynomial generators given by the exterior powers \(\Lambda^i(V)\) and \(\Lambda^i(V,\phi)\), respectively, where \(V\) is the \(2n\)-dimensional standard representation.
\begin{lemma}\label{lem:R-Sp-alternative-generators}
  In \(\K_0(\Rep(\Sp_{2n}))\), \(\Lambda^i V = E_{\omega_i} + \sum_{x\colon x < \omega_i} n_xE_x\) for all \(i\in\{1,\dots,n\}\), where the sum runs over a finite number of dominant weights \(x\) such that \(x < \omega_i\) in the partial order of \Cref{def:serres-partial-order}, and \(n_x\in \N\). In particular,
  \(
    \K_0(\Rep(\Sp_{2n})) \cong \Z[V,\Lambda^2V, \dots, \Lambda^n V]
  \).
\end{lemma}
\begin{proof}
  We first note that the fundamental weights \(\omega_i\) of \(\Sp_{2n}\) satisfy
  \begin{equation}\label{eq:Sp-fundamental-weights-order}
    \omega_{i+2} > \omega_i
  \end{equation}
  for all \(i\), i.e.\ for all \(i \in \{1,\dots,n-2\}\).  Indeed, we see from \cite[Planche~III]{bourbakiLie456} that \(\omega_{i+2}-\omega_i = \alpha_{i+1} + (\sum_{j=i+2}^{r-1}\alpha_j) + \alpha_r\), where \(\alpha_1\), \dots, \(\alpha_n\) are the simple roots, so \(\omega_{i+2}-\omega_i > 0 \).

  In characteristic zero, \(\Lambda^iV\) decomposes as
  \begin{equation}\label{eq:Sp-V-decomposition-0}
    \Lambda^iV \cong E_{\omega_i} \oplus E_{\omega_{i-2}} \oplus \dots \oplus
      \begin{cases}
        E_{\omega_1} & \text{ if \(i\) is odd } \\
        k          & \text{ if \(i\) is even }
      \end{cases}
  \end{equation}
  for \(i=1,\dots,n\) \emph{as a representation}, not just as an element of \(\K_0(\Rep(\Sp_{2n}))\).
  (See \cite[Chapter~VIII, §\,13.3 (IV)]{bourbakiLie789}
  for a proof of the corresponding statement for the Lie algebra of \(\Sp_{2n}\), and recall that in characteristic zero we can pass from decompositions of a representation of a split semisimple simply connected algebraic group to decompositions of the associated representation of its Lie algebra and vice versa using \cite[\parasign\,3.2]{Humphreys:LAG}.)
  The claim of the lemma is therefore immediate from \eqref{eq:Sp-fundamental-weights-order} in characteristic zero.

  Suppose now that \(k\) is a field of positive characteristic.
  Note that \(\Sp_{2n}\) lifts to a split reductive group scheme \(\Sp_{2n,\Z}\) over \(\Z\) \cite[Th\'eorème~7.2.0.46]{CalmesFasel:SGA3}, hence the results of \cite[section 3.7]{serre68} apply.
  Let us temporarily write \(\Sp_{2n,\Q}\) and \(\Sp_{2n,k}\) for the the symplectic groups over \(\Q\) and \(k\), which can be obtained from \(\Sp_{2n,\Z}\) via base change,
  and let us denote by \(E_{\omega_i,\Q}\) and \(E_{\omega_i,k}\) the respective irreducible representations of highest weight \(\omega_i\).
  Consider the category \(\Rep(\Sp_{2n,\Z})\) of \(\Sp_{2n,\Z}\)-representations that are finitely generated and free over \(\Z\), and the associated K-ring \(\K_0(\Rep \Sp_{2n,\Z})\), denoted 
\(\mathrm{R}_{\Z}(\Sp_{2n})\) in \cite[\parasign\,3.3]{serre68}.

  The base-change homomorphisms 
  \begin{align*}
    i_\Q\colon \K_0(\Sp_{2n,\Z}) & \to \K_0(\Sp_{2n,\Q})\\
    i_k\colon \K_0(\Sp_{2n,\Z}) & \to \K_0(\Sp_{2n,k})
  \end{align*}
 are homomorphisms of \(\lambda\)-rings, essentially because exterior powers are compatible with base change \cite[Proposition~A2.2\,(b)]{eisenbud:ca}.
  By \cite[Théorème~5]{serre68}, \(i_\Q\) is an isomorphism.  Better still, we see from \cite[Lemme~2\,(a)]{serre68} that for each representation \(E_\Q\) of \(\Sp_{2n,\Q}\) there exists a representation \(E_\Z\), not just a virtual representation, such that \(i_\Q(E_\Z) = E_\Q\). In particular, for each dominant weight \(\omega\), we have a representation \(E_{\omega,\Z}\) of \(\Sp_{2n,\Z}\) with \(i_\Q(E_{\omega,\Z}) = E_{\omega,\Q}\).
  Note, however, that there is no reason to assume that \(i_k(E_{\omega,\Z}) = E_{\omega,k}\) in general.  Rather, arguing as in \cite[\parasign\,3.8, Remarque~3]{serre68}, we find that 
  \begin{equation}\label{eq:base-change-of-irreducbile-rep}
    i_k(E_{\omega,\Z}) = E_{\omega,k} + \textstyle\sum_{x < \omega} n_x E_{x,k},
  \end{equation}
  where the sum is over all dominant \(x\) with \(x < \omega\), and only finitely many of the coefficients \(n_x \in\N_0\) are non-zero. Indeed, let us write \(\ch_k\) for the character homomorphism \(\ch_{\Sp_{2n,k}}\) (c.f. \Cref{lemma:highestweightofproductissumofweights}).  As we have a maximal torus of \(\Sp_{2n}\) defined over \(\Z\), and as restriction to this maximal torus is compatible with base change, we find that \(\ch_{\Sp_{2n,k}}(i_k(E_{\omega,\Z})) = \ch_{\Sp_{2n,\Q}}(E_{\omega,\Q})\).  Thus, \cite[Lemma~5\,(a)]{serre68}, applied both over \(\Q\) and over \(k\), yields the following two identities, for certain coefficients \(a_x, b_x \in\N_0\):
  \begin{align*}
    \ch_k(i_k(E_{\omega,\Z})) &= e^\omega + \textstyle\sum_{x < \omega}a_x e^x\\
    \ch_k(E_{\omega,k}) &= e^\omega + \textstyle\sum_{x < \omega} b_x e^x
  \end{align*}
  \Cref{eq:base-change-of-irreducbile-rep} follows from a comparison of these two identities and the injectivity of \(\ch_k\).
  
  On the other hand, the standard representation of \(\Sp_{2n}\) is already defined over \(\Z\), so we have a representation \(V_\Z\) of \(\Sp_{2n,\Z}\) with \(i_k(V_\Z) = V_k\) for arbitrary fields~\(k\).  
  Equation~\eqref{eq:Sp-V-decomposition-0}, appropriately decorated with subscripts \((-)_\Z\), is therefore equally valid in \(\K_0(\Sp_{2n,\Z})\).
  The claimed decomposition of \(\Lambda^i(V)\) in \(\K_0(\Rep(\Sp_{2n}))\) now follows by combining \eqref{eq:Sp-fundamental-weights-order}, this integral version of \eqref{eq:Sp-V-decomposition-0}, and \eqref{eq:base-change-of-irreducbile-rep}.

  The final claim is immediate from this decomposition and \Cref{rem:representationringispolynomialring-alternative-generators}.
\end{proof}

\begin{lemma}\label{lem:GW-Sp-alternative-generators}
  Let \((V,\phi)\) denote the \(2n\)-dimensional standard representation of \(\Sp_{2n}\) with its canonical anti-symmetric form.
  The exterior powers \(\Lambda^i(V,\phi)\) for \(i\in\{1,\dots,n\}\) define homogeneous polynomial generators of the \(\pm\)-graded \(\GW^\pm(k)\)-algebra \(\GW^\pm(\Rep \Sp_{2n})\):
  \[
    \GW^\pm(\Rep(\Sp_{2n})) = \GW^\pm(k)[(V,\phi),\Lambda^2(V,\phi), \dots, \Lambda^n(V,\phi)]
  \]
  The generator \(\Lambda^i(V,\phi)\) is of degree \((-1)^i\).
\end{lemma}
\begin{proof}
  Applying \Cref{lem:symmetric-rep-decomposition} to the decomposition of \(\Lambda^i(V)\) given in \Cref{lem:R-Sp-alternative-generators}, we find that 
  \(
    \Lambda^i(V,\phi) = (E_{\omega_i},\phi_{\omega_i}) + (\text{smaller terms})
  \)
  for certain symmetric isomorphism \(\phi_{\omega_i}\) on \(E_{\omega_i}\).
  So the claim is immediate from \Cref{rem:gwringispolynomialring-alternative-generators}.
 \end{proof}

\subsubsection*{Restriction to diagonal}
The group  \(\Sp_{2n}\) has a canonical subgroup isomorphic to \(\Sp_2^{\times n} = \Sp_2 \times \dots \times \Sp_n\).
This is most easily seen by replacing the anti-symmetric form \(J\) from \eqref{eq:model-of-Sp} by the isometric form \(nH^-=H^-\perp\dots\perp H^-\); then \(\Sp_2^{\times n}\) is simply given by \(n\) diagonally concatenated copies of $\Sp_2$ inside \(\Sp_{2n}\).
The inclusion of this subgroup induces a commutative diagram
\begin{equation*}
    \begin{tikzcd}
        \K_0(\Rep(\Sp_{2n})) \arrow[r, "\res"] \arrow[dr, hook]
            & \K_0(\Rep(\Sp_2^{\times n})) \arrow[d, hook] \\
        {}
            & \K_0(\Rep(T)),
    \end{tikzcd}
\end{equation*}
which shows that the restriction map $\res: \K_0(\Rep(\Sp_{2n})) \ra \K_0(\Rep(\Sp_2^{\times n}))$ is injective.  Its image can easily be identified, as follows (cf.\ explicit calculations in \cite[Appendix D]{karpenko21}). Consider the symmetric group \(S_n\) acting by permutation on \(\Sp_2^{\times n}\), and its induced action on \(\K_0(\Rep(\Sp_2^{\times n}))\).  The image of the restriction map is precisely the fixed ring of this action:
\begin{equation*}
    \K_0(\Rep(\Sp_{2n})) \xrightarrow[\res]{\cong} \K_0(\Rep(\Sp_2^{\times n}))^{S_n}
\end{equation*}
We will now show that the corresponding statement for \(\GW^\pm\) also holds.
\begin{lemma}
  The restriction map \(\GW^\pm(\Rep(\Sp_{2n})) \to \GW^\pm(\Rep(\Sp_2^{\times n}))\) is injective, with image the invariant subring under the permutation action of the symmetric group \(S_n\).
\end{lemma}
\begin{proof}
  Note that $\Sp_2^{\times n}$ is a simply connected reductive group, and that \Cref{prop:gwringispolynomialring} applies just as it applies to \(\Sp_{2n}\).
  So \(\GW^\pm(\Rep(\Sp_2^{\times n}))\) is a polynomial ring in $n$ variables, which we can identify with \(\prod_{i=1}^n\GW^\pm(\Rep(\Sp_2))\).
  Writing \(v^{(i)} := (V^{(i)},\phi^{(i)})\) for the standard representation of the \(i^{\text{th}}\) factor in \(\Sp_2^{\times n}\), equipped with its canonical anti-symmetric form, we thus obtain:
  \begin{equation}\label{eq:GW-of-diagonal-Sp2s}
    \GW^\pm(\Rep(\Sp_2^{\times n})) \cong  \Z[v^{(1)},\dots,v^{(n)}].
  \end{equation}
  Under this isomorphism, the permutation action of \(S_n\) on the left corresponds to the obvious permutation action on the generators \(v^{(i)}\) on the right.

  Now let \((V,\phi)\) denote the \(2n\)-dimensional standard representation of \(\Sp_{2n}\), equipped with its canonical anti-symmetric form.
  As we have seen in \Cref{lem:GW-Sp-alternative-generators}
  we can take the exterior powers \(\Lambda^k(V,\phi)\) for \(k = 1,\dots, n\) as polynomial generators of \(\GW^\pm(\Rep(\Sp_{2n}))\).
  Write \(\sigma^k(v^\bullet)\in\GW^\pm(\Rep(\Sp_2^{\times n}))\) for the \(k^{\text{th}}\) elementary symmetric function in the classes \(v^{(i)} = (V^{(i)},\phi^{(i)})\) introduced above.
  We claim that
  \begin{equation}\label{eq:78RW20}
    \res\Lambda^k(V,\phi) = \sigma_k(v^\bullet) + \text{ a polynomial in \(\sigma_j(v^\bullet)\) with \(j<k\)}
  \end{equation}
  To this end, recall that the operations \(\lambda^k\) defined in terms of the exterior powers \(\Lambda^k\) give \(\GW^\pm(\Rep(G))\) the structure of a pre-\(\lambda\)-ring (for any affine group scheme \(G\)), see \cite{zibrowius15} and \cite[Proposition~4.2.1]{fasel23stable}.
  For \(k=1\), equation \eqref{eq:78RW20} simply follows from the fact that \((V,\phi)\) restricts to the direct sum of the representations \((V^{(i)},\phi^{(i)})\):
  \begin{align*}
    \res(V,\phi)
    = \textstyle\sum_{i=0}^n v{(i)}
    = \sigma_1(v^\bullet)
  \end{align*}
  For higher \(k\), consider the power series expansion \(\lambda_t(x) := \sum_k \lambda^k(x)t^k\).
  For the two-dimensional standard representation \((V_2,\phi_2)\) of \(\Sp_2\), we find \(\Lambda^2(V_2,\phi_2) = (k,\det(\phi)) = (k,\id)\), and \(\Lambda^k(V_2,\phi_2) = 0\) for \(k>2\), so
  \[
    \lambda_t(V_2,\phi_2) = 1 + (V_2,\phi_2)t + t^2.
  \]
  As the restriction commutes with the \(\lambda\)-operations, this implies:
  \[
    \res\lambda_t(V,\phi)
    = \lambda_t(\textstyle\sum_i v^{(i)})
    = \textstyle\prod_{i=1}^n \lambda_t v^{(i)}
    = \textstyle\prod_{i=1}^n (1 + v^{(i)}t + t^2)     
  \]
  Now \(\res\lambda^k(V,\phi)\) is the coefficient of \(t^k\) in the above power series.  As the power series is invariant under the permutation action of \(S_n\), so is each coefficient.  So \(\res\lambda^k(V,\phi)\) is a polynomial in the symmetric polynomials \(\sigma^i(v^\bullet)\).  Moreover, the highest-degree monomials in the \(V^{(i)}\)'s that occur in the coeffient of \(t^k\) are precisely the monomials that occur \(\sigma^k(v^\bullet)\), and they occur with multiplicity one.  This proves \eqref{eq:78RW20}, and hence the Lemma.
\end{proof}

For the completion at the agumentation ideal, it is again convenient to reformulate the above description in terms of the “first Borel classes” defined as in \eqref{eq:GW-of-Sp2} above.  So let us write \(b^{(i)} := (V^{(i)},\phi^{(i)}) - H^-\) for the first Borel class of the \(i^{\text{th}}\) copy of \(\Sp_2\), and \(\sigma_i(b^\bullet)\) for the \(i^{\text{th}}\) elementary symmetric polynomial in these classes.  Clearly, isomorphism~\ref{eq:GW-of-diagonal-Sp2s} from the proof above can be rewritten as 
\[ \label{eq:gw-of-rep-sp2-on-borel-classes}
  \GW^\pm(\Rep(\Sp_2^{\times n})) \cong  \Z[b^{(1)},\dots,b^{(n)}].
\]
\begin{corollary} \label{corollary:symplectic-splitting-for-gw-on-rep}
  The restriction induces an isomorphism between \(\GW^\pm(\Rep(\Sp_{2n}))\) and the subring \(\GW^\pm[\sigma_1(b^\bullet),\dots,\sigma_n(b^\bullet)]\) of \(\GW^\pm(\Sp_2^{\times n})\).  Under this isomorphism, the graded augmentation ideal \(IO_{\Sp_{2n}}^\pm\) (see \Cref{def:graded-augmentation-ideal}) corresponds to the ideal generated by \(\sigma_1(b^\bullet), \dots, \sigma_n(b^\bullet)\).  In particular,
  \[
    \GW^\pm(\Rep(\Sp_{2n}))^{\wedge}_{IO_{\Sp_{2n}}^\pm} \cong \GW^\pm(k)\llbracket\sigma_1(b^\bullet), \dots, \sigma_n(b^\bullet)\rrbracket.
  \]
\end{corollary}

\section{Grothendieck-Witt rings of some classifying spaces}
\label{section:gw-rings-of-classifying-spaces}

The main purpose of this chapter is to prove Theorem \ref{theorem:asc-for-gw-of-sp}.
First, we introduce the specific model category we will be working with to model motivic spaces. 
Next, in \Cref{subsection:acceptable-gadgets}, we recall the notion of \emph{acceptable gadgets} as introduced in \cite[Definition 8.3]{panin18} and establish some important properties of these.
In \Cref{subsection:hermitian-asc-symplectic-groups}, we finally use the results of \Cref{subsection:some-details-symplectic-groups} to construct the Atiyah-Segal completion map, and to prove the main theorem.
In the final two subsections, we discuss the construction of a classifying space for the multiplicative group with non-trivial involution, and conjectural tools that may be used to generalize our main results to base schemes with non-trivial action and to arbitrary closed subgroups of~$\Sp_{2n}$.

\subsection{Čech localization}
    \label{subsection:cech-localization}

Let $S$ be a scheme of finite type over a field $k$. 
Our standing assumption that $\cha(k) \neq 2$ is not necessary for any of the results of this section.
Let $\sPre(\Sm_S)$ be the model category of simplicial presheaves on $\Sm_S$ with the global injective model structure, and $L_{\mot}\sPre(\Sm_S)$ its motivic localization, which is defined as $L_{\mot} = (L_{\A^1}L_{\Nis})^{\infty}$, where the repeated localizations are necessary to ensure that $L_{\mot}F$ is both Nisnevich local and $\A^1$-local for any simplicial presheaf $F$. 
This is a presheaf variant of the model category constructed in \cite{morel99}, with homotopy category the unstable motivic homotopy category $\mc{H}(S)$.  
The weak equivalences in $L_{\mot}\sPre(\Sm_S)$ will be called \emph{motivic weak equivalences}; these are precisely the maps that become isomorphisms in $\mc{H}(S)$. 

More generally, for a subcanonical topology $\tau$ on $\Sm_S$, let $L_{\tau}\sPre(\Sm_S)$ be the (left) Bousfield localization with respect to covering sieves for \(\tau\), as considered in \cite[Section~3.1]{asok17}.
This localization is also considered in \cite[Appendix~A]{dugger04}, where it is called the Čech localization because it is a localization with respect to Čech covers in the topology $\tau$, see \cite[Theorem~A5]{dugger04}.
This is the naming convention we will follow.

The fibrant objects of $L_{\tau}\sPre(\Sm_S)$ are those fibrant simplicial presheaves in $\sPre(\Sm_S)$ that satisfy $\tau$-descent; we will refer to these as \(\tau\)-fibrant. 
The weak equivalences in $\sPre(\Sm_S)$ will be called \emph{objectwise weak equivalences} and the weak equivalences in $L_\tau\sPre(\Sm_S)$ will be called \emph{$\tau$-local weak equivalences}. 
We denote by $L_{\tau}$ the $\tau$-fibrant replacement functor, viewed as endofunctor on $\sPre(\Sm_S)$.

\begin{remark}[Relation to model structure used by Morel \& Voevodsky]\label{rem:relation-to-MV-model-structure}
  Note that $L_{\tau}\sPre(\Sm_S)$ is \emph{not} the same as the Jardine model structure \(L_{\hyptau}\sPre(\Sm_S)\),
  which is Quillen equivalent to the Joyal model structure on simplicial sheaves used in \cite{morel99}.   In general, we only have succesive Bousfield localizations
\begin{equation*}
    \sPre(\Sm_S) \lra L_{\tau}\sPre(\Sm_S) \lra L_{\hyptau}\sPre(\Sm_S).
\end{equation*}
However, the Bousfield localization \(L_{\tau}\sPre(\Sm_S) \lra L_{\hyptau}\sPre(\Sm_S)\) is a Quillen equivalence for \(\tau\) the Nisnevich topology and \(S\) any reasonable base scheme.
See \cite[Remark~3.1.4]{asok17} for more details.
\end{remark}

\begin{remark}[\(\infty\)-language]\label{rem:infty-language}
We may also view $\sPre(\Sm_S)$ as an $\infty$-category as in e.g.\ \cite{hoyois17six}.
The fibrant replacement functor $L_{\tau}$ is analogous to the localization endofunctor defined in \cite[Proposition~3.4]{hoyois17six} and has the same formal properties.
Lemma~2.1 from \cite{hoyois20cdh} will play a key role in our discussion of classifying spaces (see the proof of \Cref{lemma:nice-gadgets-etale-classifying-space} below).
\end{remark}

We will need very few specifics about the \(\tau\)-local model structure.
The proof of the following proposition relies on the two subsequent lemmas, both presumably well-known.
\begin{proposition}\label{prop:schemes-and-filtered-colims-are-fibrant}\mbox{}~
  \begin{itemize}
  \item[(a)]
    For any subcanonical topology \(\tau\), every \(S\)-scheme is \(\tau\)-fibrant.
  \item[(b)]
    For any topology $\tau$ at least as coarse as the fppf topology, any filtered colimit in \(\sPre(\Sm_S)\) of \(S\)-schemes is \(\tau\)-fibrant.
  \end{itemize}
  In particular, for any scheme \(X\) as in (a) or any filtered colimit \(X\) as in (b), the \(\tau\)-fibrant replacement map $X \ra L_{\tau}(X)$ is an objectwise weak equivalence.
\end{proposition}

\begin{proof}
  For (a), let \(X\) be an \(S\)-scheme.  We need to show that \(X\) is fibrant in the injective model structure on \(\sPre(\Sm_S)\) and \(\tau\)-local.  For the first assertion, see \Cref{lemma:constant-sset-lifts-for-weak-equivalences} below.  For the second assertion, see \cite[Example~3.1.2]{asok17}: a simplicially constant presheaf is \(\tau\)-local if and only if it is a \(\tau\)-sheaf.  The presheaf \(X\) in question is a \(\tau\)-sheaf precisely because \(\tau\) is assumed to be subcanonical. 
  For the filtered colimit in (b), the claim follows by the same argument plus the observation, spelled out in \Cref{lem:colim-etale} below, that such filtered colimits are indeed sheaves.
The final assertion is just an application of Ken Brown's Lemma.
\end{proof}

As remarked above, the following two lemmas are probably well-known to the experts.  \Cref{lemma:constant-sset-lifts-for-weak-equivalences} is stated without proof in \cite[section 2.6]{rezk10}. We include full details here for future reference.

\begin{lemma} \label{lemma:constant-sset-lifts-for-weak-equivalences}
Let $\mc{C}$ be a small category.
A simplicially constant  presheaf $F \in \sPre(\mc{C})$ is fibrant in any model structure on $\sPre(\mc{C})$ whose weak equivalences are objectwise weak equivalences.
In particular, constant simplicial presheaves are fibrant in the injective model structure.
\end{lemma}

\begin{proof}
Let $i: A \ra B$ be an objectwise weak equivalence in $\sPre(\mc{C})$ and $a: A \ra F$ a map of simplicial presheaves.
To prove the statement, it suffices to prove the existence of a map $b: B \ra F$ such that $bi = a$.
Note that $F$ is an objectwise Kan complex and therefore projective fibrant, so that $a$ factors through a functorial projective fibrant replacement $RA$ of $A$, as in the following diagram
\begin{equation*}
    \begin{tikzcd}
        A \arrow[r, "\sim"] \arrow[d, "i"] & RA \arrow[r, "a'"] \arrow[d, "Ri"] & F \\
        B \arrow[r, "\sim"] & RB \arrow[ur, dashed, swap, "b'"],
    \end{tikzcd}
\end{equation*}
and a lift $b$ of $a$ along $i$ exists if a lift $b': RB \ra F$ of $a'$ along $Ri$ exists.
Therefore, we may and do assume both $A$ and $B$ to be projective fibrant, that is, they are both objectwise Kan complexes.

We first produce a collection of maps $b_C:B(C) \ra F(C)$ such that $a_C = b_C i_C$ for each $C \in \mc{C}$.
Let $C \in \mc{C}$.
Since $i_C: A(C) \ra B(C)$ is a weak equivalence of Kan complexes, it is a simplicial homotopy equivalence.
Let $j_C: B(C) \ra A(C)$ be a homotopy inverse to $i_C$ and define $b_C: B(C) \ra F(C)$ as the composition $a_C j_C$. 
Then the composition $a_C j_C i_C: A(C) \ra F(C)$ is homotopic to the map $a_C$, and since homotopic maps into a constant simplicial set are equal, $a_C = b_C i_C$.
Thus, for each $C \in \mc{C}$ we have $b_C = a_C j_C$ for a fixed choice of homotopy inverse $j_C$ of $i_C$.

We now need to show that the square
\begin{equation*}
    \begin{tikzcd}
        B(D) \arrow[r, "b_D"] \arrow[d, "B(f)"] & F(D) \arrow[d, "F(f)"] \\
        B(C) \arrow[r, "b_C"] & F(C)
    \end{tikzcd}
\end{equation*}
commutes for any arrow $f: C \ra D$ in $\mc{C}$. 
Consider the diagram
\begin{equation*}
    \begin{tikzcd}
        B(D) \arrow[r, "j_D"] \arrow[d, "B(f)"] & A(D) \arrow[d, "A(f)"] \arrow[r, "a_D"] & F(D) \arrow[d, "F(f)"] \\
        B(C) \arrow[r, "j_C"] & A(C) \arrow[r, "a_C"] & F(C),
    \end{tikzcd}
\end{equation*}
and note that the right square commutes.
Again because homotopic maps into a constant simplicial set are equal, it suffices to show that the left square commutes up to homotopy.
Note that
\begin{equation*}
    j_C i_C A(f) j_D = j_C B(f) i_D j_D
\end{equation*}
because $i$ is a map of simplicial presheaves.
Moreover, the left hand side is homotopic to $A(f) j_D$ and the right hand side is homotopic to $j_C B(f)$ because $j_C i_C \sim \id_{A(C)}$ and $i_D j_D \sim \id_{B(D)}$ by definition.
Hence, $F(f) b_D = b_C B(f)$, as was to be shown.
\end{proof}

\begin{lemma}\label{lem:colim-etale}
Consider a topology $\tau$ coarser than the fppf topology.  Given a family $(F_i)_{i \in I}$ of representable $\tau$-sheaves on a quasi-compact scheme \(S\), where $I$ is a filtered index category, the filtered colimit $F_{\infty}:= \colim_{i \in I}F_i$ is also a $\tau$-sheaf.
\end{lemma} 

\begin{proof}
  As $S$ is quasi-compact, we have to check the sheaf condition for $F_{\infty}$ only for finite $\tau$-covers $\U$ of a scheme $U$, by the following argument.
  Let $\U$ be an arbitrary $\tau$-cover of $U$.
  As $U$ is quasi-compact, there exists a finite Zariski cover $\mc{V}$ of $U$.
  Note that $\U \cap \mc{V}$ is a refinement of $\U$.
  For $V \in \mc{V}$, we consider the $\tau$-cover $\U \cap V$ of $V$. 
  By (the proof of) \stackstag{021P}, this refines to a finite $\tau$-cover of $V$.
  Hence $\U \cap \mc{V}$ refines to a finite $\tau$-cover of $U$, 
  and from now on we assume that $\U$ is a finite $\tau$-cover of $U$.
We denote the associated covering sieve by $h_{\U}$, which is a subpresheaf of the representable presheaf $h_U$. 
  We are reduced to showing that if the canonical map $\Hom(h_U,F_i) \to \Hom(h_{\U},F_i)$ is a bijection for all $i \in I$, then it also is a bijection for $i=\infty$. 
  The representable sheaf $h_U$ is compact, hence commutes with filtered colimits. 
  Moreover, for a finite covering $\U=(U_j \to U)_{j \in J}$, $h_{\U}$ can be described by the explicit small coequalizer 
  \begin{equation*}
      \coprod_{j,j'} U_j \times _U U_{j'} \rightrightarrows \coprod_j U_j
  \end{equation*}
  of representables, hence is also compact, and the claim follows. 
\end{proof}

\subsection{Acceptable gadgets and classifying spaces}
    \label{subsection:acceptable-gadgets}

Let $S$ be a scheme of finite type over a field $k$. 
Our standing assumption that $\cha(k) \neq 2$ is not necessary for any of the results of this section.

\paragraph{Classifying spaces}
We now fix a $\tau$-sheaf of groups \(G\) over $S$.  
In applications, $G$ will often be a linear algebraic group over $\Spec k$ base changed to $S$.

The \emph{simplicial classifying space} \(\B{G}\) is defined in \cite{morel99} as the nerve of \(G\) viewed as a presheaf of groupoids with a single object (so the $n$-simplices are $\B{G}_n := G^{\times n}$ and face and degeneracy maps are given by composition and inserting identities, as usual).

\begin{definition}\label{def:BetG}
  For a topology \(\tau\) at least as coarse as the fppf topology, the \emph{\(\tau\)-classifying space} $\B_{\tau}G$ of $G$ is defined as $\B_{\tau}G := L_{\tau}(\B{G})$.
\end{definition}

Note that \(\B_{\tau}G\) is well-defined -- irrespective of any particular choice of fibrant replacement functor -- up to objectwise weak equivalence.
We now briefly discuss the dependency of \(\B_{\tau} G\) on the topology \(\tau\).
For the definition of $\tau$-locally trivial $G$-torsors, we refer to \cite[Definition~2.3.1]{asok18}.
\begin{lemma}\label{BG-up-to-objectwise-we}
  The following conditions on \(G\) and \(\tau\) are equivalent:
  \begin{enumerate}
  \item[1.] Every (fppf-locally trivial) \(G\)-torsor in \(\Sm_S\) is already \(\tau\)-locally trivial.
  \item[2.] The canonical map $\B_{\tau} G \to \B_{\fppf}G$ is an objectwise weak equivalence.
  \end{enumerate}
\end{lemma}
\begin{proof}
  The stated equivalence follows from \cite[Lemma~2.3.2\,(i) and (ii)]{asok18}.
  (See also \cite[Lemma 4.1.18]{morel99}).
\end{proof}

\begin{example}[smooth affine groups]
  For smooth affine group schemes $G$, every $G$-torsor is already a étale-locally trivial, so  $\B_{\fppf}G$ is objectwise weakly equivalent to $\B_{\et}G$.
  Indeed, smoothness is preserved by faithfully flat descent, so any $G$-torsor is smooth over the base, and smooth morphisms admit sections étale locally.
\end{example}

\begin{example}[special groups]\label{rem:special-groups}
  A linear algebraic group \(G\) over \(S\) is called \emph{special} if every étale-locally trivial \(G\)-torsor over a (not necessarily smooth) \(S\)-scheme is already Zariski-locally trivial.  
  Thus, for a special group \(G\), \(\B_\Zar G\) is objectwise weakly equivalent to \(\B_\et G\), and to \(\B_\tau G\) for any intermediate topology such as \(\tau = \Nis\).
  Special groups include, in particular, 
  split tori,
  \(\GL_r\),
  \(\SL_{2r}\),
  \(\Sp_{2r}\),
  and finite products of these \cite[Lemmas 3.1 and 3.2]{reichsteintossici20}.
  (In \cite{serre58,reichsteintossici20}, the notion is introduced and discussed only for groups defined over a field \(k\), but note that the defining property is stable under base change along \(S \to \Spec(k)\).  Special groups over fields have been fully classified \cite{grothendieck58,huruguen16,merkurjev22}.)
\end{example}
Of course, for \(\tau = \Nis\) or any coarser topology, \(\B_{\tau}G\) is \emph{motivically} equivalent to the simplicial classifying space \(\B G\), for any \(G\), as \(\tau\)-local equivalences are motivic equivalences in this case. Only for finer choices of topologies \(\tau\), \(\B_\tau G\) may be motivically distinct. 

\begin{remark}[Comparison with \cite{morel99}]
  By construction, \(\B_{\tau}G\) only satisfies $\tau$-descent. However, it also satisfies $\tau$-hyperdescent for dimension reasons by the argument given in the proof of \cite[Lemma~2.3.2]{asok18}, so by Ken Brown's Lemma, the canonical map $\B_{\tau}G \ra L_{\hyptau}\BG$ is an objectwise weak equivalence.  In view of \Cref{rem:relation-to-MV-model-structure}, this shows, in particular, that the definition of $\B_{\et}G$ used here agrees with the definition in \cite[p.130]{morel99}.
\end{remark}

\paragraph{Outline}
Using the \emph{admissible gadgets with a nice $G$-action} of \cite[Definition~4.2.1, Definition~4.2.4]{morel99}, one constructs a geometric approximation of a universal \(G\)-torsor over \(\B_{\et}G\),
namely a colimit of $G$-torsors in \(\Sm_S\).
We will instead use the more flexible concept of \emph{acceptable gadgets}, following \cite[Definition 8.3]{panin18}, and adding a version of a nice $G$-action in this setting.
Suppose that we have a commutative diagram of $G$-torsors
\begin{equation}\label{diag:torsor-sequence}
    \begin{tikzcd}
        X_1 \arrow[r] \arrow[d]
            & X_2 \arrow[r] \arrow[d]
                & \dots \\
        Y_1 \arrow[r]
            & Y_2 \arrow[r]
                & \dots
    \end{tikzcd}
\end{equation}
in $\Sm_S$, where the horizontal arrows are closed immersions, such that the colimit of this diagram should be thought of a motivic approximation of $EG \ra \B_{\et}G$, with $EG$ being a motivically contractible space with free $G$-action. We will make this idea precise further in a more general setting for any Grothendieck topology $\tau$ at least as coarse as the fppf topology,
and for $\tau$-locally trivial $G$-torsors and $\tau$-classifying spaces $\B_{\tau}G$.

Namely, we will establish a chain of weak equivalences
\[
  Y_\infty \underset{(1)}\simeq L_{\tau} (\squotient{X_\infty}{G}) \underset{(2)}\simeq L_{\tau} (\squotient{S}{G}) \underset{(3)}\cong B_{\tau} G 
\]
All quotients \(\squotient{-}{G}\) appearing here are stacky; see \Cref{def:stacky-quotient}.  The objectwise weak equivalence~(1) is established below in \Cref{lemma:colimit-cofibs-of-torsors-is-torsor}.
Equivalence~(2) is a motivic weak equivalence, see \Cref{lemma:nice-gadgets-etale-classifying-space}. 
The isomorphism~(3) is completely general and clear from the definitions, see \Cref{BG-as-stacky-quotient}. 
Combined, these equivalences yield \emph{a motivic weak equivalence \(Y_\infty \simeq \B_{\tau}G\) for any sequence of \(\tau\)-locally trivial \(G\)-torsors as in \eqref{diag:torsor-sequence} in which the sequence of \(X_i\)'s forms an acceptable gadget.}
This result, which is stated as Theorem \ref{generalBGtheorem} below, will be used in the next section to show that the products $(\HP^n)^{\times r}$ approximate $\BSp_{2}^{\times r}$, and similarly certain quaternionic Grassmannians approximate $\BSp_{2r}$.
These are crucial steps in the proof of Atiyah-Segal completion for $\BSp_{2r}$.

We begin by recalling the definition of an acceptable gadget.

\begin{definition} \label{definition:acceptable-gadget}
An \emph{acceptable gadget} $(X_n)_{n \in I}$ over a scheme $S$ is a countable totally ordered set $I$ together with a diagram $X: I \ra \Sm_S$ such that each $X_n = X(n)$ is quasi-projective and for each $n < m$ in $I$, the map $X(n < m): X_n \ra X_{m}$ is a closed immersion of $S$-schemes, and such that for any Henselian regular local ring $R$ and any commutative square
\begin{equation*}
    \begin{tikzcd}
        \partial \Delta^i_R \arrow[d] \arrow[r] 
            & X_n \arrow[d] \\
        \Delta^i_R \arrow[r]
            & S,
    \end{tikzcd}
\end{equation*}
there exists an $m \geq n$ and a map $\Delta^i_R \ra X_m$ making the diagram
\begin{equation} \label{eq:gadget-lift}
    \begin{tikzcd}
        \partial \Delta^i_R \arrow[d] \arrow[r] 
            & X_n \arrow[r]
                & X_m \arrow[d] \\
        \Delta^i_R \arrow[rr] \arrow[urr, dashed]
            & { } 
                & S
    \end{tikzcd}
\end{equation}
commute.
\end{definition}

Here is a lemma that shows that acceptable gadgets yield contractible motivic spaces, which will make them useful to construct contractible spaces with free actions, whose quotients are classifying spaces as we will see in \Cref{lemma:nice-gadgets-etale-classifying-space}.
Some version of \Cref{lemma:colim-acceptable-gadget-is-contractible} is implied in \cite[p. 954 after Proposition 8.5]{panin18}.
\begin{lemma} \label{lemma:colim-acceptable-gadget-is-contractible}
Let $(X_i)_{i \in I}$ be an acceptable gadget over a scheme $S$. 
Then $X_{\infty} = \colim_i X_i$ is contractible in the category $\mc{H}(S)$ of motivic spaces over $S$.
\end{lemma}

\begin{proof}
  By \cite[Lemma~2.3.8]{morel99} and \cite[Lemma~3.1.11]{morel99}, it suffices to show that the simplicial set $\Sing(X_{\infty})(R)$ is contractible for every regular Henselian local ring \(R\) over \(S\), in other words, for every Nisnevich point of $S$, see also \Cref{rem:relation-to-MV-model-structure}.
  
By construction, a map $\partial \Delta^n \ra \Sing(X_{\infty})(R)$ is given by a morphism $\partial \Delta^n_R \ra X_{\infty}$.
Let $f$ be such a morphism.
Since $\partial \Delta^n_R$ is representable by an object of $\Sm_S$, it is compact, and $f$ factors through a finite stage $f': \partial \Delta^n_R \ra X_i$ for some $i \in I$. 
Since $(X_i)_{i \in I}$ is an acceptable gadget, there exists $j > i$ in $I$ such that the composition $\partial \Delta^n_R \ra X_i \ra X_j$ extends to a morphism $\Delta^n_R \ra X_j$, which shows that $\Sing(X_{\infty})(R)$ is contractible, and the proof is finished.
\end{proof}

We prove a lemma that allows us to construct acceptable gadgets from existing ones.
The first of these shows that acceptable gadgets are stable under base change.

\begin{lemma}\label{propertiesofacceptablegadgets}
  Let $(X_n)_{n \in \N}$ and $(Y_n)_{n \in \N}$ an acceptable gadgets over a scheme $S$. 
  \begin{itemize}
  \item[(a)]  \label{lemma:acceptable-gadget-base-change}
    For any $Y \in \Sm_S$, the base change $(X_n \times_S Y)_{n \in \N}$ is an acceptable gadget over $Y$.  \item[(b)] \label{lemma:acceptable-gadget-cofinal}
    For any cofinal $I \subset \N$, the sequence $(X_n)_{n \in I}$ is an acceptable gadget over $S$.
  \item[(c)]\label{lemma:product-gadget-is-acceptable}
    The product $(X_n \times_S Y_n)_{n \in \N}$ is an acceptable gadget over $S$.
  \end{itemize}
\end{lemma}

\begin{proof}
We only prove (c); the proofs of (a) and (b) are similar.  
  Consider a commutative diagram
\begin{equation*}
    \begin{tikzcd}
        \partial \Delta^i_R \arrow[d] \arrow[r] 
            & X_n \times_S Y_n \arrow[d] \\
        \Delta^i_R \arrow[r]
            & S.
    \end{tikzcd}
\end{equation*}
The projection maps $X_n \times_S Y_n \ra X_n$ and $X_n \times_S Y_n \ra Y_n$ yield commutative diagrams
\begin{equation*}
    \begin{tikzcd}
        \partial \Delta^i_R \arrow[d] \arrow[r] 
            & X_n \arrow[r]
                & X_{m_1} \arrow[d] \\
        \Delta^i_R \arrow[rr] \arrow[urr, dashed]
            & { } 
                & S
    \end{tikzcd}
    \qquad
    \begin{tikzcd}
        \partial \Delta^i_R \arrow[d] \arrow[r] 
            & Y_n \arrow[r]
                & Y_{m_2} \arrow[d] \\
        \Delta^i_R \arrow[rr] \arrow[urr, dashed]
            & { } 
                & S,
    \end{tikzcd}
\end{equation*}
as $(X_n)_{n \in \N}$ and $(Y_n)_{n \in \N}$ are acceptable gadgets.
Let $m = \max\{m_1,m_2\}$.
Then the compositions $\Delta^i_R \ra X_{m_1} \ra X_m$ and $\Delta^i_R \ra Y_{m_2} \ra Y_m$ yield a map $\Delta^i_R \ra X_m \times_S Y_m$ making the diagram
\begin{equation*}
    \begin{tikzcd}
        \partial \Delta^i_R \arrow[d] \arrow[r] 
            & X_n \times_S Y_n \arrow[r]
                & X_m \times_S Y_m \arrow[d] \\
        \Delta^i_R \arrow[rr] \arrow[urr, dashed]
            & { } 
                & S
    \end{tikzcd}
\end{equation*}
commute, as was to be shown.
\end{proof}

  All quotients \(\squotient{-}{G}\) used in the following will be stacky quotients, defined as follows.  Recall that the action groupoid of a set \(X\) with a group action is the groupoid with objects the elements of \(X\), and morphisms from \(x_1\in X\) to \(x_2\in X\) given by the group elements that send \(x_1\) to \(x_2\).
\begin{definition}\label{def:stacky-quotient}
  For a presheaf of sets \(X\) on \(\Sm_S\) with \(G\)-action, the \emph{stacky quotient} \(\squotient{X}{G}\) is the simplicial presheaf given by the nerve of the action groupoid of~\(X\).
\end{definition}

We will mainly be interested in two extremal cases:  the stacky quotient of the (trivial) \(G\)-action on the base scheme~\(S\), and the stacky quotient of the (free) \(G\)-action on a \(G\)-torsor.  These are discussed in the following two lemmas, respectively.

\begin{lemma}\label{BG-as-stacky-quotient}
  For any subcanonical topology \(\tau\), we have a canonical isomorphism of simplicial presheaves
  $\B_\tau G \cong L_{\tau}(\squotient{S}{G})$.
\end{lemma}

\begin{proof}
      For any \(U\in \Sm_S\), the action groupoid of the \(G(U)\)-action on the one-point set \(S(U)\) is isomorphic to the groupoid \(G(U)\) used in the definition of \(B G\).
\end{proof}

\begin{lemma}\label{lem:torsor-as-stacky-quotient}
  For any $\tau$-locally trivial $G$-torsor $\pi\colon X \ra Y$ in $\Sm_S$, we have a canonical objectwise weak equivalence $L_{\tau}(\squotient{X}{G}) \ra Y$.
\end{lemma}
\begin{proof}
  \newcommand{\prequotient}[2]{#1/_{\!\mathrm{pre}\,}#2}
  For any presheaf of sets \(X\) with a free \(G\)-action, the stacky quotient $\squotient{X}{G}$ is objectwise weakly equivalent to the presheaf quotient $U \mapsto X(U)/G(U)$, which we temporarily denote by \(\prequotient{X}{G}\).
  Indeed, when the \(G(U)\)-action on \(X(U)\) is free, the action groupoid \(X(U)\) is canonically equivalent to the orbit set \(X(U)/G(U)\), viewed as a discrete groupoid, so we obtain a weak equivalence when passing to nerves.
  
    For a smooth $S$-scheme $X$ with a free $G$-action, we hence have an objectwise weak equivalence from $L_{\tau}(\squotient{X}{G})$ to the $\tau$-sheafification \((\prequotient{X}{G})^\tau\) of the presheaf quotient.
    Indeed, note first that the objectwise weak equivalence \(\squotient{X}{G}\to\prequotient{X}{G}\) described above induces an objectwise weak equivalence \(L_\tau(\squotient{X}{G})\to L_\tau(\prequotient{X}{G})\) (by two-out-of-three and Ken Brown's Lemma). 
    Secondly, as \((\prequotient{X}{G})^\tau\) is \(\tau\)-fibrant (by \Cref{lemma:constant-sset-lifts-for-weak-equivalences} and \cite[Example~3.1.2]{asok17}), the canonical map \(\prequotient{X}{G}\to (\prequotient{X}{G})^\tau\) factorizes through a morphism \(L_\tau(\prequotient{X}{G})\to (\prequotient{X}{G})^\tau\), which also is an objectwise weak equivalence by \cite[Example~3.1.2]{asok17}.
    So altogether we obtain the objectwise weak equivalence claimed above. Finally, in the situation at hand, we have an isomorphism \((\prequotient{X}{G})^\tau\cong Y\).
    This completes the proof.
  \end{proof}

We now provide a variant of \cite[Definition 4.2.4]{morel99}.

\begin{definition} \label{def:nice-acceptable-gadget}
An \emph{acceptable $G$-gadget $(X_n)_{n \in I}$ over $S$} is an acceptable gadget over $S$ satisfying:
\begin{enumerate}[label=(\roman*)]
\item for each $n \in I$, $X_n$ is endowed with a
  $G$-action; and
\item for each $n < m$ in $I$, the corresponding closed immersion $X_n \ra X_m$ is $G$-equivariant.
\end{enumerate}
\end{definition}
The following result generalizes \Cref{lem:torsor-as-stacky-quotient} from a single torsor to a sequences of torsors.
\begin{lemma} \label{lemma:colimit-cofibs-of-torsors-is-torsor}
  Let $(X_n \ra Y_n)_{n \in \N}$ be a sequence of $G$-torsors as in \eqref{diag:torsor-sequence}, with all torsors  $\tau$-locally trivial and all maps \(X_n \to X_{n+1}\) monomorphisms. Then the induced map 
  \begin{equation*}
    L_\tau(\squotient{X_\infty}{G}) \lra Y_\infty
  \end{equation*}
  is an objectwise weak equivalence.
  In particular, this holds if $(X_n)_{n \in \N}$ is an acceptable $G$-gadget with a $\tau$-locally trivial $G$-torsor $X_n \ra Y_n$ for each $n \in \N$.
\end{lemma}
\begin{proof}
  All objects are cofibrant, and the maps $X_n \ra X_m$, $\squotient{X_n}{G} \ra \squotient{X_m}{G}$ and \(Y_n\to Y_m\) are monomorphisms, hence cofibrations in our model structure.
  So \(X_\infty\), \(\colim(\squotient{X_n}{G})\) and \(Y_\infty\) are colimits of cofibrations and thus homotopy colimits. 
  The canonical maps $\squotient{X_n}{G} \ra Y_n$, which are \(\tau\)-local weak equivalences by \Cref{lem:torsor-as-stacky-quotient}, therefore induce a \(\tau\)-local weak equivalence 
  \(\colim_n (\squotient{X_n}{G}) \lra Y_{\infty}\).
  Finally, we observe that the action groupoid of $X_\infty$ is the strict colimit of the action groupoids of the $X_n$, so that we have an isomorphism
  \begin{equation*}
    \colim_n (\squotient{X_n}{G}) \cong \squotient{X_\infty}{G}.
  \end{equation*}
  Hence, we obtain a \(\tau\)-local weak equivalence \(\squotient{X_\infty}{G} \lra Y_{\infty}\).  
  As $Y_\infty$ is $\tau$-fibrant by \Cref{prop:schemes-and-filtered-colims-are-fibrant}\,(b), this \(\tau\)-local weak equivalence factors through $L_\tau(\squotient{X_\infty}{G})$,
  and by Ken Brown's Lemma we thus obtain the desired objectwise weak equivalence
  \begin{equation*}
    L_\tau(\squotient{X_\infty}{G}) \lra Y_\infty.
  \end{equation*}
  If $(X_n)_{n \in \N}$ is an acceptable $G$-gadget with a $\tau$-locally trivial $G$-torsor $X_n \ra Y_n$ for each $n \in \N$, the maps $X_n \ra X_m$ are closed immersions and therefore monomorphisms by assumption.
  As \(X_n(U) \to X_m(U)\) is a \(G(U)\)-equivariant monomorphism for every \(U\in\Sm_X\), it follows that the induced maps \(X_n(U)/G(U)\to X_m(U)/G(U)\) are also monomorphisms.
  We argue at the end of the proof of \Cref{lem:torsor-as-stacky-quotient} that \(Y_n\) and \(Y_m\) differ from the presheaf quotients only by sheafification.  As sheafification is exact, this shows that \(Y_n \to Y_m\) is also a monomorphism. 
\end{proof}

We now verify a crucial generalized contractibility condition for acceptable \(G\)-gadgets. Given a $\tau$-locally trivial torsor $T \ra Y$ in $\Sm_S$, consider the map $L_{\tau}(X_\infty \times_S T) \ra Y$ given by applying first $(\squotient{-}{G})$ and then $L_{\tau}$ to the $G$-equivariant projection map $X_\infty \times_S T \ra T$ and then composing with the map 
$L_{\tau}(\squotient{T}{G}) \ra Y$ of \Cref{lem:torsor-as-stacky-quotient}. 
This map is canonical up to the choice of \(L_{\tau}\).

\begin{proposition} \label{lemma:nice-gadget-torsor-property-for-nis-special-groups}%\label{lemma:nice-gadget-torsor-prop-new}
    Let $(X_n)_{n \in I}$ be an acceptable \(G\)-gadget over $S$.
    Then, for each \(\tau\)-locally trivial \(G\)-torsor $T \ra Y$ in $\Sm_S$, the map 
    \begin{equation*}
        L_{\tau}(\squotient{(X_{\infty} \times_S T)}{G}) \ra Y
    \end{equation*} is a motivic weak equivalence.    
\end{proposition}
\begin{proof}
    Apply \Cref{lemma:colimit-cofibs-of-torsors-is-torsor} to the sequence of $\tau$-locally trivial $G$-torsors $X_n \times_S T \ra X_n \times_S Y$ to obtain an objectwise weak equivalence
    \begin{equation*}
      L_\tau(\squotient{(X_\infty \times_S T)}{G}) \to X_\infty \times_S Y.  
    \end{equation*}
    This objectwise weak equivalence is, in particular, a motivic weak equivalence.
    As the map in the statement of the proposition can be obtained by composing this map with the projection map $p\colon X_\infty \times_S Y \ra Y$, it suffices to show that \(p\) is a motivic equivalence over~\(S\).

    We  know that \(p\) is a motivic weak equivalence over $Y$ by \Cref{lemma:colim-acceptable-gadget-is-contractible,propertiesofacceptablegadgets}.
    We now use the functor $u_{\sharp}\colon L_{\mot}\sPre(\Sm_Y) \ra L_{\mot}\sPre(\Sm_S)$ from \cite[Proposition 4.5.4]{ayoub07} for the structure map $u\colon Y \to S$ to deduce that \(p\) is also a motivic weak equivalence over~\(S\). 
    (These general results do not use the assumption on the ``coefficients'' discussed at the beginning of section 4.4 of loc.\ cit.
    Compare also \cite[section 4.1]{hoyois17six} for a related discussion in the $\infty$-setting.) 
    More precisely, $u_\sharp$ maps a smooth scheme over $Y$ to the same smooth scheme over $S$ by the proof of loc.\ cit., and is a left Quillen functor for the projective model structures by \cite[Theorem 4.5.14]{ayoub07}.
    As left adjoints preserve colimits, $u_\sharp$ also maps the morphism $p\colon X_{\infty} \times_S Y \to Y$ over $Y$ to the same morphism considered over~\(S\).
    Hence, as $u_{\sharp}$ is left Quillen, \(p\) is a motivic weak equivalence over $S$ for the projective motivic model structure, and thus also for the injective model structure, which has the same weak equivalences.
  \end{proof}

\begin{remark}\label{generalizedtorsor}
  There are definitions of $G$-torsors in the literature which are more general than \cite[Definition~2.3.1]{asok18}, see e.g.\ \cite[p.~128]{morel99},
  but which have in common that the conclusion of \Cref{lem:torsor-as-stacky-quotient} holds.  The proofs of \Cref{lemma:colimit-cofibs-of-torsors-is-torsor} and \Cref{lemma:nice-gadget-torsor-property-for-nis-special-groups} immediately generalize to these more general torsors provided the bases $Y_n$ are representable.
\end{remark}

\begin{proposition} \label{lemma:nice-gadgets-etale-classifying-space}
 For any acceptable $G$-gadget $(X_n)_{n \in I}$ over $S$, 
 $L_{\tau}(\squotient{X_{\infty}}{G})$ is motivically equivalent to $L_\tau(\squotient{S}{G})$.
\end{proposition}
\begin{proof}
  By \cite[Lemma~2.1]{hoyois20cdh} 
  it suffices to verify the condition established in \Cref{lemma:nice-gadget-torsor-property-for-nis-special-groups}.
  Note that, though  \cite[Lemma~2.1]{hoyois20cdh} is only stated for the fppf-topology, it holds equally for any finer topology~\(\tau\), as is evident from its (very short) proof.
  We also note that, while \cite[Lemma~2.1]{hoyois20cdh} is stated in terms of general \(\infty\)-\(G\)-torsors, for a simplicically discrete group \(G\) any such torsor over a simplicially discrete base $X$ is again discrete, as it is locally isomorphic to \(X\times_S G\).
  So in our setting it suffices to verify the assumptions for simplicially discrete \(G\)-torsors as in \cite[Definition~2.3.1]{asok18}.
  (Alternatively, we could use \Cref{generalizedtorsor}.)
  Finally, the $\infty$-quotients appearing in \cite[Lemma 2.1]{hoyois20cdh} coincide with our stacky quotients. See \Cref{BG-as-stacky-quotient,lem:torsor-as-stacky-quotient} above for special cases, and \cite[section~3]{nss14} for details in the general case.
\end{proof}
The following theorem is the promised precise version of the ``outline'' at the beginning of this section.

\begin{theorem}\label{generalBGtheorem}
  Let $(X_n)_{n \in \N}$ be an acceptable $G$-gadget such that we have a \(\tau\)-locally trivial \(G\)-torsor $X_n \ra Y_n$ and induced monomorphisms $Y_n \ra Y_m$ for all $n, m \in \N$. 
  Then we have motivic weak equivalences
  \begin{equation*}
    Y_{\infty} \simeq \B_{\tau}G \simeq \B_{\fppf} G. 
  \end{equation*} 
\end{theorem}
\begin{proof}
As already explained in the outline, \Cref{lemma:colimit-cofibs-of-torsors-is-torsor}, \Cref{lemma:nice-gadgets-etale-classifying-space} and \Cref{BG-as-stacky-quotient} yield a chain of motivic weak equivalences $Y_\infty \simeq L_{\tau} (\squotient{X_\infty}{G}) \simeq L_{\tau} (\squotient{S}{G}) \simeq B_{\tau} G$.  As the same chain of equivalences also applies for any topology finer than \(\tau\), we also obtain the motivic equivalence with \(\B_{\fppf} G\).
\end{proof}

As recalled in \Cref{rem:special-groups}, for \emph{special} affine algebraic group schemes \(G\) any torsor is Zariski-locally trivial, so that \Cref{generalBGtheorem} will yield motivic equivalences  \(Y_{\infty} \simeq \B G \simeq \B_\fppf{G}\).

\subsection{Hermitian ASC for symplectic groups}
    \label{subsection:hermitian-asc-symplectic-groups}

Let $G$ be a linear algebraic group over a field $k$ of characteristic not two.
Note that, in particular, $G$ is smooth and an étale sheaf on $\Sm_k$. 

\begin{definition}\label{def:graded-augmentation-ideal}
  We define \(IO_G\) and \(IO_G^\pm\) as kernels of restriction maps:
  \begin{align*}
    IO_G&:=\ker(\GW^+(\Rep(G))\to \GW^+(k))\\
    IO_G^\pm&:=\ker(\GW^\pm(\Rep(G))\to \GW^\pm(k))
  \end{align*}
\end{definition}
Our definition of \(IO_G\) agrees with the definition given in \cite{rohrbach22completion}.  The graded ideal \(IO_G^\pm\) is not considered there.
\begin{lemma}\label{lemma:graded-completion-versus-ungraded-completion}
  The \(IO^\pm_G\)-adic topology on \(\GW^\pm(\Rep(G))\) agrees with the \(IO_G\)-adic topology.
\end{lemma}
\begin{proof}
  This is a general fact about graded ideals.  Consider a \(\pm\)-graded ring \(R= R^+ \oplus R^-\) and a graded ideal \(\ideal a = \ideal a^+ \oplus \ideal a^-\subset R\).
  Clearly \((\ideal a^+)^i\cdot R \subset \ideal a^i\).  On the other hand, as \((\ideal a^-)^2 \subset \ideal a^+\), we find that:
  \[
    \ideal a^{2i} = (\ideal a^+ + \ideal a^-)^{2i} \subseteq \sum_{j=0}^{2i} \ideal (a^+)^{2i-j}\ideal (a^-)^j \subseteq (\ideal a^+)^i \cdot R^+ \oplus (\ideal a^+)^i \cdot R^-.
  \]
  So \(\ideal a^{2i} \subset (\ideal a^+)^i\cdot R\).  
  This shows that the \(\ideal a\)-adic topology on \(R\) agrees with the \(\ideal a^+\)-adic topology.
\end{proof}
Recall from \cite{panin22quaternionic,panin18} the $\Sp_{2r}$-torsors $\HU(r,n)$ over the quaternionic Grassmannians $\HGr(r,n)$ associated with the tautological symplectic bundle, defined over a smooth quasi-projective base scheme $S$.

\begin{proposition} \label{proposition-quaternionic-gadgets}
  The sequence of $\Sp_{2r}$-torsors $\HU(r,n) \to \HGr(r,n)$ with $n \in \N$ of \cite[Proposition 8.5]{panin18} defines an acceptable $\Sp_{2r}$-gadget $(\HU(r,n))_{n \in \N}$ over $S$.
\end{proposition}
\begin{proof}
  It is shown in \cite[Proposition 8.5]{panin18} that $(\HU(r,n))_{n \in \N}$ is an acceptable gadget over $S$. It obviously also satisfies the conditions of an acceptable $\Sp_{2r}$-gadget, where the compatible actions come from the structure as $\Sp_{2r}$-torsors.
\end{proof}

\begin{theorem}[Panin-Walter]\label{theorem:panin-walter-bsp}
  There are motivic weak equivalences
  $$\B_{\Nis}\Sp_{2r} \simeq \B_{\et}\Sp_{2r} \simeq \HGr(r,\infty).$$
\end{theorem}

\begin{proof}
  Both equivalences are stated in \cite[after Proposition 8.5]{panin18}, along with a brief indication on how to modify the arguments of \cite{morel99} to obtain a proof, using their concept of an acceptable gadget.
  Section~\ref{subsection:acceptable-gadgets} provides more details for this argument.
  The first equivalence is evident from the fact that \(\Sp_{2r}\) is special (see \Cref{rem:special-groups}).
  The second is immediate from \Cref{generalBGtheorem} applied to the acceptable gadget of \Cref{proposition-quaternionic-gadgets}.
\end{proof}

\begin{remark}
The motivic equivalence $\B_{\et}\Sp_{2r} \simeq \HGr(r,\infty)$ of \Cref{theorem:panin-walter-bsp} was proven in \cite[Proposition~5]{schlichting15} for $r = \infty$ using a different technique.
However, \cite[Proposition~3]{schlichting15}, which is an ingredient of this alternative proof, has no obvious analogue for finite $r$. 
\end{remark}

\begin{remark}\label{rem:comparison-with-mv}
  For the reader's convenience, let us compare the discussion above to some arguments in \cite{morel99}. \Cref{theorem:panin-walter-bsp} is an analogue for \(\Sp_{2r}\) of \cite[Proposition 4.3.7]{morel99} for $\GL_r$.  
  In either case, the first equivalence is established by noting that \(\GL_r\) and \(\Sp_{2r}\), respectively, are special, and by using \cite[Proposition 4.1.18]{morel99}. 
  The second equivalence rests on \cite[Proposition~4.2.6]{morel99}, which has its parallel in our \Cref{lemma:nice-gadgets-etale-classifying-space}.
  Both of these intermediate results rely on the contractibility of colimits of gadgets.  For admissible gadgets, this is proved in \cite[Proposition 4.2.3]{morel99}, using the ambient vector bundles, while for acceptable gadgets this is \Cref{lemma:colim-acceptable-gadget-is-contractible} above.
  In \cite{morel99}, this contractibility enters via \cite[Lemma~4.2.9]{morel99}; see the proof of \cite[Proposition~4.2.6]{morel99} spelled out below \cite[Lemma~4.2.9]{morel99}. 
  In our case, the contractibility enters via the corresponding \Cref{lemma:nice-gadget-torsor-property-for-nis-special-groups}.
  Despite all these parallels, the definition of acceptable \(G\)-gadget employed here is significantly simpler
  than the definition of nice admissible \(G\)-gadget used by Morel and Voevodsky: there is no analogue of part (iii) of \cite[Definition~4.2.4]{morel99} in our \Cref{def:nice-acceptable-gadget}.
  The main reason why this simplification is possible is \cite[Lemma~2.1]{hoyois20cdh}.
\end{remark}

For \(\Sp_2\), \Cref{proposition-quaternionic-gadgets} yields an acceptable gadget $(\HU(1,n+1))_{n \geq 1}$ built from \(\Sp_2\)-torsors $\HU(1,n+1) \to \HP^{n}$, where \(\HP^n\) is the quaternionic projective space defined in \cite{panin22quaternionic}.  In combination with \Cref{propertiesofacceptablegadgets}\,(c), we moreover obtain acceptable gadgets \((\HU(1,n+1)^{\times r})_{n \geq 1}\) built from \(\Sp_{2}^{\times r}\)-torsors over products of quaternionic projective spaces \((\HP^n)^{\times r}\), for any \(r\in \N\).

\begin{proposition}\label{GWBSp-approximated-by-GWHP}
The following maps are isomorphisms:
  \begin{enumerate}[label=(\roman*)]
      \item the canonical map \(\GW^{\pm}(\B \Sp_2) \to \lim_n \GW^{\pm}(\HP^n)\);
      \item the canonical map \(\GW^{\pm}(\B \Sp_2^{\times r}) \to \lim_n \GW^{\pm}((\HP^n)^{\times r})\); and
      \item the canonical map \(\GW^{\pm}(\B \Sp_{2r}) \to \lim_n \GW^{\pm}(\HGr(r,n))\).
  \end{enumerate}
\end{proposition}
\begin{proof}
  By \cite[Theorem 13.4]{panin18}, bigraded  $\GW$-theory is represented by a commutative ring spectrum $\BO$, which is $(8,4)$-periodic by \cite[Theorem 7.5]{panin18}. 
Restricting $\BO$ to bidegrees $(0,0)$ and $(4,2)$ yields the commutative ring $\GW^{\pm}(X)$ for motivic spaces $X$ which are not schemes.
  For an acceptable $G$-gadget $(X_n)_{n \in \N}$ such that we have compatible $G$-torsors $X_n \ra Y_n$ for each $n$, we may consider the motivic Milnor exact sequences (see e.g.\ \cite[Theorem~5.7]{panin22}):
\begin{align*}
0 \to \textstyle\lim^1_n \BO^{-1,0}(Y_n) \to \GW^{+}(Y_{\infty}) \to \lim_n \GW^{+}(Y_n) \to 0 & \\
0 \to \textstyle\lim^1_n \BO^{3,2}(Y_n) \to \GW^{-}(Y_{\infty}) \to \lim_n \GW^{-}(Y_n) \to 0 &
\end{align*}
These yield a canonical isomorphisms
\(
    \GW^{\pm}(Y_{\infty}) \to \lim_n \GW^{\pm}(Y_n)
\)
provided the $\lim^1$-terms vanish.
When \(Y_n = \HGr(r,n)\), they do vanish, as \cite[Theorem~9.5]{panin18}
provides the necessary surjections of $\BO$-groups, thus proving (i) and (iii).
To prove (ii), we argue similarly, using \cite[Theorem~9.5]{panin18} inductively by viewing $(\HP^n)^{\times r}$ as a trivial quaternionic projective bundle over $(\HP^n)^{\times r - 1}$. 
\end{proof}

In view of this \nameCref{GWBSp-approximated-by-GWHP}, the reader preferring motivic spaces to actual varieties may replace several $\lim_n \GW^{\pm}(\HP^n)$ below by $\GW^{\pm}(\B_{\et}\Sp_2)$ or $\GW^{\pm}(\BSp_2)$.
\medskip

We are now ready to examine the case $G = \Sp_2$ in detail.
The gadget $(\HU(1,n+1))_{n \geq 1}$ has nice formal properties for \emph{Borel classes}.
Let $(\mc{U}_n, \psi_n)$ be the tautological rank $2$ symplectic bundle on $\HP^n$.
Note that it defines an element in $\GW^-(\HP^n)$.
We define the universal first Borel class 
\begin{equation*}
    b_1 = (b_{1,(n)})_{n \in \N} \in \lim_n \GW^{-}(\HP^n)
\end{equation*}
to be the limit of the first Borel classes $b_{1,(n)} \in \GW^{-}(\HP^n)$,
where $b_{1,(n)} = b_1(\mc{U}_n, \psi_n)$ is the first Borel class of $(\mc{U}_n, \psi_n)$ on $\HP^n$, as defined in \cite[Definition~8.3]{panin22quaternionic}. 
By \cite[Proposition~9.9]{panin18}, we have 
\begin{equation*}
  b_{1,(n)} = [\mc{U}_n, \psi_n] - H^- \in \GW^{-}(\HP^n),
\end{equation*}
where $H^-$ is the symplectic hyperbolic space of rank~$2$ with trivial $\Sp_2$-action.
Let $(V, \phi)$ be a trivial symplectic vector space of rank $2n + 2$ in which $(V_0, \phi_0)$ is a distinguished rank $2$ symplectic subspace, and write $\Sp_{2n + 2} = \Sp(V, \phi)$.
In this case $(V, \phi)$ is the standard representation of $\Sp_{2n + 2}$ discussed in \Cref{expl:sp}. 
There is a decomposition
\begin{equation*}
    (V, \phi) \cong (V_0, \phi_0) \perp (V_0^{\perp}, \phi_0^{\perp}).
\end{equation*}
Then $\Sp_{2n + 2}$ acts on $\HP^n = \HP(V, \phi)$, and the stabilizer of the distinguished point $(V_0, \phi_0)$ in $\HP^n$ is given by those $a \in \Sp(V, \phi)$ that fix $(V_0, \phi_0)$ and $(V_0^{\perp}, \phi_0^{\perp})$, and is therefore given by $\Sp_{2n} \times \Sp_2$, where we make the identifications $\Sp_{2n} = \Sp(V_0^{\perp}, \phi_0^{\perp})$ and $\Sp_2 = \Sp(V_0, \phi_0)$.
Under the identification
\begin{equation*}
    \HP^n \cong \Sp_{2n+2}/(\Sp_{2n} \times \Sp_2),
\end{equation*}
the canonical \(\Sp_2\)-torsor \(\HU(1,n+1)\) is given by \(\Sp_{2n+2}/\Sp_{2n}\),
which corresponds to the universal bundle $(\mc{U}_n, \psi_n)$ under the correspondence between $\Sp_2$-torsors and rank $2$ symplectic bundles as outlined in e.g.\ \cite[Section~3.3, p.~1025]{asok18}.
We let $X_n$ denote $\HU(1,n+1)$ for ease of notation.
Let $\pi_n: X_n \ra \Spec k$ be the $G$-equivariant structure map.
Since the pullback functor $\gamma_n^*: \Vect(\HP^n) \ra \Vect^G(X_n)$ along the projection $\gamma_n: X_n \ra \HP^n$ is an equivalence of exact categories with duality, $\gamma_n^*$ induces a canonical isomorphism $\gamma_n^*: \GW^{\pm}(\HP^n) \ra \GW^{\pm}_G(X_n)$.
The following lemma studies the Atiyah-Segal map $\GW^{\pm}(\Rep(\Sp_2)) \ra \GW^{\pm}(\BSp_2)$ induced by approximations $\pi_n: X_n \to \Spec k$ of the pullback along $\ESp_2 \ra \Spec k$.

\begin{proposition} \label{proposition:borel-class-standard-representation} Let $H^{-}$ be the trivial symplectic plane bundle equipped with the trivial $\Sp_2$-action and $[V_0,\phi_0]$ the $\Sp_2$-representation from above. 
The composition 
\begin{equation*}
    (\gamma_n^*)^{-1}\pi_n^*: \GW^{\pm}(\Rep(\Sp_2)) \lra \GW^{\pm}_{\Sp_2}(X_n) \lra \GW^{\pm}(\HP^n)   
\end{equation*}
sends $[V_0,\phi_0] - H^{-}$ to the (first) Borel class $b_{1,(n)}= [\mc{U}_n, \psi_n] - H^- $.
\end{proposition}

\begin{proof}
The pullback map $\pi_n^*: \GW^{\pm}(\Rep(\Sp_2)) \ra \GW^{\pm}_{\Sp_2}(X_n)$ sends $[V_0,\phi_0]$ to the class $[V_n, \phi_n]$ of the trivial rank $2$ symplectic bundle on $X_n$ with the standard $\Sp_2$-action.
Since $(\mc{U}_n, \psi_n)$ is a rank $2$ symplectic bundle corresponding to the torsor $X_n = \Symp_{(\mc{U}_n,\psi_n)}$, its pullback to $X_n$ is also the trivial symplectic bundle $[V_n, \phi_n]$ with the standard $\Sp_2$-action.
As $X_n \ra \HP^n$ is an $\Sp_2$-torsor, the map
\begin{equation*}
    \Sp_2 \times_k X_n \ra X_n \times_{\HP^n} X_n
\end{equation*}
given on points by $(a,x) \mapsto (ax, x)$ is an isomorphism.
Hence, $\Aut_{\Symp}(H_{X_n}^-) \cong \Sp_2 \times_k X_n \cong \gamma_n^*\Symp_{(\mc{U}_n,\psi_n)} \cong \Symp_{(f^*\mc{U}_n,f^*\psi_n)}$ and we deduce
\begin{equation*}
    \pi_n^*([V_0,\phi_0] - H^-) = [V_n,\phi_n] - H^- = \gamma_n^*([\mc{U}_n, \psi_n] - H^- ),
\end{equation*}
as was to be shown.
\end{proof}

As an immediate corollary, we obtain Atiyah-Segal completion for $G=\Sp_2$ and $\GW^{\pm}$. 

\begin{corollary} \label{corollary:asc-for-sp2}
For $G=\Sp_2$, the map $\GW^{\pm}(\Rep(G)) \ra \lim_n \GW^{\pm}(\HP^n)$ from \Cref{proposition:borel-class-standard-representation} above is a completion of $\GW^{\pm}(\Rep(G))$ with respect to the Hermitian augmentation ideal $IO_G$. 
\end{corollary}
\begin{proof} The computation (\ref{eq:GW-of-Sp2}) immediately implies that $IO^{\pm}_{\Sp_2}$ is generated by $[V_0,\phi_0] - H^{-}$. 
By \cite[Theorem 9.5]{panin18},
\begin{equation*}
    \lim_n \GW^{\pm}(\HP^n) \cong \GW^{\pm}(k)[[b_1]].
\end{equation*}
Hence the claim follows from \Cref{proposition:borel-class-standard-representation} and \Cref{lemma:graded-completion-versus-ungraded-completion}.
\end{proof}

Now we prove Atiyah-Segal completion for general $\Sp_{2r}$.
We will use the diagram
\begin{equation} \label{diagram:symplectic-splitting}
    \begin{tikzcd}
        \GW^{\pm}(\Rep(\Sp_{2r})) \arrow[r, "\res"]
            & \GW^{\pm}(\Rep(\Sp_2^{\times r})) \arrow[d] \\
        \GW^{\pm}(\BSp_{2r}) \arrow[r]
            & \GW^{\pm}(\BSp_2^{\times r}),
    \end{tikzcd}
\end{equation}
in which the vertical arrow is a completion by \Cref{prop:gw-rep-to-gw-bsp-is-completion}.
We start with the following result of Panin and Walter.

\begin{proposition} \label{proposition:gw-of-product-of-quaternionic-proj}
For $1 \leq i \leq r$, let $y_i \in \lim_n \GW^{\pm}((\HP^n)^{\times r})$ be the element defined by the inverse limit 
\begin{equation*}
    lim_{n \in \N}b_1(\mc{U}_n^{(i)}, \phi_n^{(i)})
\end{equation*}
of the first Borel classes 
$b_1(\mc{U}_n^{(i)}, \phi_n^{(i)})$
of the $i$-th tautological rank $2$ bundle on $(\HP^n)^{\times r}$.
Then
\begin{equation*}
    \lim_n \GW^{\pm}((\HP^n)^{\times r}) = \GW^{\pm}(k)[[y_1, \dots, y_r]].
\end{equation*}
\end{proposition}

\begin{proof}
For $r = 1$, this is a consequence of \cite[section 11]{panin22quaternionic} or \cite[Theorem 9.5]{panin18} as already recalled above.
Note that $(\HP^n)^{\times r} \ra (\HP^n)^{\times r - 1}$ is a trivial $\HP^n$-bundle, so by \cite[Theorem 9.4]{panin18},
\begin{equation*}
    \GW^{\pm}((\HP^n)^{\times r}) \cong \frac{\GW^{\pm}((\HP^n)^{\times r - 1})[b_1(\mc{U}_n^{(r)}, \phi_n^{(r)})]}{(b_1(\mc{U}_n^{(r)}, \phi_n^{(r)})^n)}.
\end{equation*}
Iterating, we obtain
\begin{equation*}
    \GW^{\pm}((\HP^n)^{\times r}) \cong \frac{\GW^{\pm}(k)[b_1(\mc{U}_n^{(1)}, \phi_n^{(1)}), \dots, b_1(\mc{U}_n^{(r)}, \phi_n^{(r)})]}{(b_1(\mc{U}_n^{(1)}, \phi_n^{(1)})^n, \dots, b_1(\mc{U}_n^{(r)}, \phi_n^{(r)})^n)}
\end{equation*}
and it follows that
\begin{equation*}
    \lim_n \GW^{\pm}((\HP^n)^{\times r}) = \GW^{\pm}(k)[[y_1, \dots, y_r]],
\end{equation*}
as was to be shown.
\end{proof}

For $n \in \N$, let $f_n: (\HP^n)^{\times r} \ra \HGr(r,rn)$ be the canonical map such that the pullback of the tautological rank $2r$ symplectic bundle on $\HGr(r, rn)$ is the orthogonal sum of the rank $2$ symplectic bundles $(\mc{U}_n^{(i)}, \phi_n^{(i)})$ on $(\HP^n)^{\times r}$; such a map exists by the universal property of quaternionic Grassmannians discussed in \cite[10-11]{panin22quaternionic}.
Note that $(\HGr(r,rn))_{n \in \N}$ is an acceptable gadget by \Cref{lemma:acceptable-gadget-cofinal}.
Recall also that we have $\lim_n \GW^{\pm}(\HGr(r,rn)) \cong \GW^{\pm}(k)[[b_1,...,b_r]]$ by \cite[Theorem 11.4]{panin22quaternionic}.
\begin{theorem} \label{theorem:symplectic-splitting-for-gw}
Let $b_i \in \lim_n \GW^{\pm}(\HGr(r,rn))$ be the element defined by the inverse limit
\begin{equation*}
    b_1(\mc{U}_{r,rn}, \phi_{r,rn})_{n \in \N}
\end{equation*}
of Borel classes of the tautological rank $2r$ bundle on $\HGr(r,rn)$.
The limit 
\begin{equation*}
    f^*: \lim_n \GW^{\pm}(\HGr(r,rn)) \lra \lim_n \GW^{\pm}((\HP^n)^{\times r})
\end{equation*}
of the pullback maps $f_n^*$ sends $b_i$ to the $i$-th symmetric elementary polynomial in the variables $y_j$ defined in \Cref{proposition:gw-of-product-of-quaternionic-proj}.
\end{theorem}

\begin{proof}
The $i$-th tautological symplectic bundle of rank $2$ on $(\HP^n)^{\times r}$ is an orthogonal direct summand of $f_n^*(\mc{U}_{r,rn}, \phi_{r,rn})$ by definition of $f_n$.
We can now run the same argument as in the first part of the proof of \cite[Theorem 10.2]{panin22quaternionic} to show that the image of $b_i$ is the $i$-th symmetric elementary polynomial in the variables $y_j$.
\end{proof}

Recall from \Cref{GWBSp-approximated-by-GWHP} that we have an isomorphism
\[
  \GW^{\pm}(\BSp_2^{\times r}) \cong \lim_n \GW^{\pm}((\HP^n)^{\times r})
\]
with the ring structure given by \Cref{proposition:gw-of-product-of-quaternionic-proj}.
Considering the composition $(\gamma_n)^{-1}\circ \pi_n^*$ and $n \to \infty$ for $G=\Sp_2^{\times r}$ and $G=\Sp_{2r}$, we obtain morphisms of $\GW^{\pm}(k)$-algebras $\GW^{\pm}(\Rep(G)) \ra \GW^{\pm}(\B G)$
as well, both generalizing (take $r=1$) the map of Corollary \ref{corollary:asc-for-sp2}.
\begin{proposition} \label{prop:gw-rep-to-gw-bsp-is-completion}
The map of $\GW^{\pm}(k)$-algebras
\begin{equation*}
    \GW^{\pm}(\Rep(\Sp_2^{\times r})) \ra \GW^{\pm}(\BSp_2^{\times r})
\end{equation*}
defined above exhibits $\GW^{\pm}(\BSp_2^{\times r})$ as the completion of $\GW^{\pm}(\Rep(\Sp_2^{\times r}))$ with respect to $IO_{\Sp_2^{\times r}}$. 
\end{proposition}

\begin{proof}
By \Cref{proposition:borel-class-standard-representation}, the generator $b^{(i)} \in \GW^{\pm}(\Rep(\Sp_2^{\times r}))$ of (\ref{eq:gw-of-rep-sp2-on-borel-classes}) is mapped to the Borel class $y_i \in \GW^{\pm}(\BSp_2^{\times r})$ of \Cref{proposition:gw-of-product-of-quaternionic-proj} for each $1 \leq i \leq r$. 
Since the ideal $IO_{\Sp_2^{\times r}} \subset \GW^{\pm}(\Rep(\Sp_2^{\times r}))$ is generated by the classes $b^{(i)}$, the result follows.
\end{proof}

We are now ready to prove the following Atiyah-Segal completion result for Hermitian $K$-theory and symplectic groups:

\begin{corollary} \label{corollary:asc-for-gw-of-sp}
The map of $\GW^{\pm}(k)$-algebras
\begin{equation*}
    \GW^{\pm}(\Rep(\Sp_{2r})) \lra \GW^{\pm}(\BSp_{2r})
\end{equation*}
defined above exhibits $\GW^{\pm}(\BSp_{2r})$ as the completion of $\GW^{\pm}(\Rep(\Sp_{2r}))$ with respect to $IO_{\Sp_{2r}}$. 
\end{corollary}

\begin{proof}
By Corollary~\ref{corollary:symplectic-splitting-for-gw-on-rep}, the map $\res: \GW^{\pm}(\Rep(\Sp_{2r})) \ra \GW^{\pm}(\Rep(\Sp_2^{\times r}))$ from the upper line of (\ref{diagram:symplectic-splitting}) is injective and maps (higher) Borel classes  to elementary symmetric polynomials in the generators $b^{(i)} \in \GW^{\pm}(\Rep(\Sp_2^{\times r}))$.
Similarly, by \Cref{theorem:symplectic-splitting-for-gw} and \Cref{GWBSp-approximated-by-GWHP}, the map $\GW^{\pm}(\BSp_{2r}) \ra \GW^{\pm}(\BSp_2^{\times r})$ from (\ref{diagram:symplectic-splitting}) is injective and maps (higher) Borel classes to elementary symmetric polynomials in the generators $y_i \in \GW^{\pm}(\BSp_2^{\times r})$.
It follows that the image of $IO_{\Sp_{2r}}$ in $\GW^{\pm}(\Rep(\Sp_2^{\times r}))$ is
\begin{equation*}
    IO_{\Sp_{2r}} = IO_{\Sp_2^{\times r}} \cap \GW^{\pm}(\Rep(\Sp_{2r}))
\end{equation*}
if we consider $\GW^{\pm}(\Rep(\Sp_{2r}))$ as a subalgebra of $\GW^{\pm}(\Rep(\Sp_2^{\times r}))$ via the restriction map.
Thus it follows from \Cref{prop:gw-rep-to-gw-bsp-is-completion} that the map
\begin{equation*}
    \GW^{\pm}(\Rep(\Sp_{2r})) \lra \GW^{\pm}(\BSp_{2r})
\end{equation*}
exhibits $\GW^{\pm}(\BSp_{2r})$ as the completion of $\GW^{\pm}(\Rep(\Sp_{2r}))$ with respect to $IO_{\Sp_{2r}}$, as was to be shown.
\end{proof}

\subsection{Classifying space for multiplicative group with a non-trivial involution}
    \label{subsection:classifying-space-group-with-involution}

Consider the multiplicative group \(\Gm\) with the involution \(\iota\colon t\mapsto t^{-1}\), as in \Cref{eg:torusinv}.
In this section, we show how to approximate the classifying space \(\B{\Gm}\) in a way that is compatible with the involution.  That is, we will construct \(\Gm\)-torsors \(U_n \to B_n\) such that the torsors \(U_n\) form an acceptable gadget, and are equipped with involutions \(\iota\colon U_n\to U_n\) that are compatible with the \(\Gm\)-action in the sense that
\begin{equation}\label{eq:Gminv-compatibility}
  \iota(t.x) = \iota(t).\iota(x)
\end{equation}
for \(t\in\Gm\), \(x\in U_n\). This condition ensures, in particular, that the involution on \(U_n\) descends to an involution on \(B_n\). 

\begin{remark}\label{rem:Gminv-usual-approximation-does-not-work}
  Usually, \(\B{\Gm}\) is approximated by the projective spaces \(\P^n\).  We have \(\Gm\)-torsors \(\A^{n+1}\setminus 0 \to \P^n\), and these torsors form an acceptable gadget.   However, there seems to be no involution \(\iota\) on \(\A^{n+1}\setminus 0\) satisfying \eqref{eq:Gminv-compatibility}.
\end{remark}
Concretely, we will use the following principal \(\Gm\)-torsors:
\[
  U_n := \{ (x,y)\in \A^{n+1}\times \A^{n+1}\mid \transpose{x}y = 1 \}
\]
with \(\Gm\)-action \(t.(x,y) := (tx, t^{-1}y)\) and involution \(\iota\colon (x,y)\mapsto (y,x)\).  This involution clearly satisfies \eqref{eq:Gminv-compatibility}.
Quotienting by the \(\Gm\)-action, we obtain the following open subschemes of \(\P^n\times \P^n\):
\[
  B_n := \{ ([x],[y])\in\P^n\times\P^n\mid \transpose{x}y \neq 0 \}
\]
The induced involution on \(B_n\) is given by \(\iota\colon ([x],[y])\mapsto ([y],[x])\).  

\begin{remark}\label{rem:Gminv-identification-with-GL-quotients}
To verify that the obvious projection \(U_n\to B_n\) is a principal \(\Gm\)-torsor, we can identify it with the canonical projection 
\[
  \GL_{1+n}/(1\times\GL_n) \to \GL_{1+n}/(\Gm\times\GL_n)
\]
via the map that sends
the left coset represented by a matrix \(A\) to the pair \((x,y)\) consisting of the first column of \(A\), and of the first row of \(A^{-1}\).
\[
  \begin{tikzcd}
    \frac{\GL_{1+n}}{1\times \GL_n} \arrow[r,"\cong"] \arrow[d,"\Gm"] & U_n \arrow[d,"\Gm"] \\
    \frac{\GL_{1+n}}{\Gm\times \GL_n} \arrow[r,"\cong"] & B_n
  \end{tikzcd}
\]
The involutions on \(U_n\) and \(B_n\) are induced by the involution \(A\mapsto \transpose{(A^{-1})}\) on \(\GL_{1+n}\).
\end{remark}
\begin{remark}
  The fixed points of \(B_n\) under the involution can be identified with the scheme \(\{[x]\in\P^n \mid \transpose{x}x \neq 0 \}\), i.e.\ with the complement of a quadric in \(\P^n\).  This complement is often used as an algebro-geometric replacement for real projective space \(\R\P^n\), for example in \cite{delzant62} or \cite{zibrowius14}.  It is a special case of the algebraic orthogonal Grassmanians of Schlichting and Tripathi \cite{schlichting15}.
\end{remark}
\begin{lemma}\label{lem:Gminv-acceptable-gadget}
  The schemes \(U_n\) with the obvious inclusions \(U_n\subset U_{n+1} \subset \dots\) form an acceptable gadget in the sense of \Cref{definition:acceptable-gadget}.
\end{lemma}
\begin{proof}
  Let \(R\) be an arbitrary commutative ring, and \(g\in R\).  As in the proof of \cite[Proposition~8.5]{panin18}, it suffices to show that, given an arbitrary morphism \(\Spec(R/g)\to U_n\), we can fill in the dashed arrow in the following diagram:
  \[
    \begin{tikzcd}
      \Spec(R/g) \arrow[d,hookrightarrow] \arrow[r]& U_n \arrow[r,hookrightarrow] & U_{n+1} \\
      \Spec(R) \arrow[urr,dashed]
    \end{tikzcd}
  \]
  The given arrow corresponds to a tuple \((\bar a_0,\dots,\bar a_n,\bar b_0,\dots,\bar b_n)\) in \((R/g)^{2n+2}\) such that \(\sum_i \bar a_i \bar b_i = 1\) in \(R/g\). Pick an arbitrary lift \((a_0,\dots, a_n, b_0,\dots,b_n)\in R^{2g+2}\). Then the previous equality tells us that there exists some element \(r\in R\) such that \(\sum_i a_i b_i = 1 + gr\) in \(R\).  The composition from \(\Spec(R/g)\) to \(U_{n+1}\) corresponds to the tuple \((\bar a_0,\dots,\bar a_n,\bar 0,\bar b_0,\dots,\bar b_n,\bar 0)\) in \((R/g)^{2n+4}\).  In order to construct the dashed arrow in a way that the diagram commutes, we need to construct a lift of this tuple, say,
  \[
    (a_0', \dots, a_n', g c, b_0',\dots,b_n', gd )
  \]
  such that \(\sum_i a_i'b_i' + g^2 cd = 1\).  Pick \(a_i' := a_i\), \(b_i' := b_i - grb_i\), and \( c := d := r\).  A quick calculation shows that this choice fits the bill.
\end{proof}

\subsection{Tools for generalizations to base schemes with non trivial group actions}
Recall that if we have Atiyah-Segal completion for $\GL_n$ for schemes $X$ with arbitrary $\GL_n$-action, then applying it to $X=\GL_n/H$ yields the completion theorem for $H$ over a point $Spec(k)$. The same applies to $\Sp_{2n}$, and hence either case would cover all $H$ which are e.g. split reductive, and in particular all finite groups.

We now briefly recall two techniques that have been succesfully used to generalize results on schemes $X$ with trivial $G$-action to general $G$-schemes $X$. Both only apply to $G=T$ a torus, and both have the same underlying idea: under suitable assumptions, there is a big open $G$-subscheme $U$ of $X$ on which the action of $X$ has a very simple product description. Using this product description, we may prove the desired result for $U$, and then using a finite number of induction steps also for $X$, assuming that the equivariant cohomology theory we care about (here: equivariant Hermitian $K$-theory) satisfies a suitable equivariant localization theorem.

The first technique is to work with $T$-filtrable schemes for $G=T$ an ``algebraic torus'', and has been used e.g.\ by Brion and Krishna, and more recently by \cite{tabuada21}.
This goes back to Bialynicki-Birula. The main ingredient is probably \cite[Theorem 2.5]{bialynicki73}, which states that $U$ is $T$-equivariantly isomorphic to $(U \cap X^T) \times V$, where $V$ is a finite-dimensional $T$-module.  

The second technique uses the ``torus generic slice theorem'' of Thomason, see \cite[Proposition 4.10]{thomason86}. Here $X$ is very general, and $G=T$ is a ``diagonalizable torus''. The main geometric result here is that we have a $T$-equivariant isomorphism $U \cong T/T' \times U/T$ with $T'$ a diagonalizable subtorus, and the induction is then done essentially on page 804 of loc. cit..

Both techniques would be useful if we could generalize them from tori $T$ to $\GL_n$ or $\Sp_{2n}$, or to products of $\Sp_2$, or to tori with involution. We unfortunately don't have this yet. Still, the above method is expected to yield 
the completion theorem for subgroups of $T$. e.g.\ products of groups of roots of unity $\mu_l$. If we assume that these roots $\mu_l$ are contained in the base field, we the could deduce the completion theorem for the corresponding direct sums of (constant) cyclic groups.

\printbibliography

\end{document}